\documentclass[12pt]{article}
\usepackage{xr}
\usepackage{amsmath,epsfig,amssymb,amsfonts,amsthm,epsfig,verbatim, enumitem}
\usepackage{color,graphicx,lscape,longtable,natbib}
\usepackage{wasysym}

\newtheorem{Th}{\underline{\bf Theorem}}

\newtheorem{Rem}{\underline{\bf Remark}}

\newtheorem{Pro}{Proposition}

\def\bse{\begin{eqnarray*}}
  \def\ese{\end{eqnarray*}}
\def\be{\begin{eqnarray}}
  \def\ee{\end{eqnarray}}
\def\bsq{\begin{equation*}}
  \def\esq{\end{equation*}}
\def\bq{\begin{equation}}
  \def\eq{\end{equation}}

\def\var{\hbox{var}}

\def\wh{\widehat}
\def\wt{\widetilde}

\def\n{\nonumber}
\def\bias{\mbox{bias}}

\def\MSE{\mbox{MSE}}

\def\sumi{\sum_{i=1}^n}

\def\trans{^{\rm T}}

\def\bb{{\boldsymbol\beta}}

\def\0{{\bf 0}}
\def\1{{\bf 1}}
\def\A{{\bf A}}

\def\e{{\bf e}}
\def\R{{\bf R}}

\def\a{{\bf a}}
\def\B{{\bf B}}

\def\b{{\bf b}}

\def\T{{\bf T}}

\def\S{{\bf S}}

\def\W{{\bf W}}

\def\T{{\bf T}}

\def\S{{\bf S}}

\def\H{{\bf H}}

\def\Ubar{{\overline{U}}}

\def\bq{\begin{equation}}
  \def\eq{\end{equation}}

\def\wh{\widehat}
\def\wt{\widetilde}
\def\trans{^{\rm T}}

\def\log{{\rm log}}
\def\bias{\hbox{bias}}

\def\squarebox#1{\hbox to #1{\hfill\vbox to #1{\vfill}}}
\def\btheta{{\boldsymbol \theta}}

\def\var{\hbox{var}}

\def\bse{\begin{eqnarray*}}
  \def\ese{\end{eqnarray*}}
\def\be{\begin{eqnarray}}
  \def\ee{\end{eqnarray}}
\def\bsq{\begin{equation*}}
  \def\esq{\end{equation*}}
\def\bq{\begin{equation}}
  \def\eq{\end{equation}}
\def\wh{\widehat}
\def\wt{\widetilde}
\def\bias{\hbox{bias}}

\def\trans{^{\rm T}}

\def\boxit#1{\vbox{\hrule\hbox{\vrule\kern6pt\vbox{\kern6pt#1\kern6pt}\kern6pt\vrule}\hrule}}

\def\boxit#1{\vbox{\hrule\hbox{\vrule\kern6pt  \vbox{\kern6pt#1\kern6pt}\kern6pt\vrule}\hrule}}

\begin{document}

%\begin{frontmatter}

\pagenumbering{arabic}
\begin{center}
{\Large\bf A spline-assisted semiparametric approach to
  nonparametric measurement error models}
\end{center}

\begin{center}
Fei Jiang\\
Department of Statistics, University of Hongkong, Hongkong\\
feijiang@hku.hk\\
{\it and}\\
Yanyuan Ma \\
Department of Statistics, Penn State University, University Park, PA
16802\\
yzm63@psu.edu\\
{\it and} \\ Raymond J. Carroll\\
Department of Statistics, 
Texas A\&M University, 
College Station TX and \\
School of Mathematical and Physical Sciences, 
University of Technology Sydney, 
Broadway NSW 2007, Australia\\
carroll@stat.tamu.edu\\
% Yang Li\\
% Department of Statistics, Renmin University, Beijing, China\\
% yang.li@ruc.edu.cn\\

\end{center}

\begin{abstract}
It is well known that the minimax rates of convergence of
nonparametric density and regression function estimation of a
random variable measured with error is much slower than the rate in the error free
case. Surprisingly, we show that  if one is willing to impose a
relatively mild assumption in
requiring that the error-prone variable has a compact support, then
the results can be greatly improved. We describe new and
constructive methods to take full advantage of the compact support
assumption via spline-assisted semiparametric methods.
We further prove that the new estimator achieves
the usual nonparametric
rate in estimating both the density and regression functions as if
there were no measurement error. The
proof involves linear and bilinear operator theories, semiparametric
theory, asymptotic analysis regarding Bsplines, as well as
integral equation treatments.
The performance of the new methods is demonstrated through
several simulations and a data example.
\end{abstract}

{\it Keywords:} Errors in covariates, measurement error,
semiparametrics, spline.

\section{Introduction}

Density estimation is a familiar problem in the nonparametric estimation
literature. Generally, we observe independent and identically
distributed variables $X_1, \dots, X_n$ from a
distribution with probability density function (pdf) $f_X(x)$ and
nonparametric estimators such as kernel methods are available in the
literature to estimate $f_X(x)$. Even when the $X_i$'s are not
directly observed, nonparametric estimation of $f_X(x)$ can still be
carried out based on their surrogates. Specifically,
assume that instead of observing $X_i$,
we observe $W_i\equiv X_i+U_i$, where $U_i$ is a mean zero random error
independent of $X_i$ and follows a
distribution with pdf $f_U(u)$. This problem has been studied extensively
in the literature
\citep{carrollhall1988,liutaylor1989,
  stefanskicarroll1990,zhang1990,fan1991} and it is well known
that the estimator of $f_X(x)$ may converge
very slowly. For example, when the error distribution $f_U(u)$ is
normal or within the class of  ``super smooth'' distributions, an
 estimator can only converge to $f_X(x)$ at the
rate of $\{\log(n)\}^{-k}$ where $k$ is a positive constant.
When the error distribution $f_U(u)$ is Laplace or another ``ordinary smooth''
type, the convergence rate is $n^{-k}$, where $k$ is a
positive constant smaller than 0.25. Here ``super smooth'' and
``ordinary smooth'' are characteristics well explained in, for
example, \cite{fan1991}. These convergence rates are minimax
for general $X$. When the $X_i$'s are observed, nonparametric density
estimation usually performs much better than
these results. We describe a simple constraint that allows rates of
convergence that are much faster than the minimax rates, and 
achieves the same rate as when there is no measurement error.

A  problem parallel to nonparametric density estimation is
nonparametric regression. Likewise, when observations $(X_i, Y_i)$,
$i=1, \dots, n$, are available, many nonparametric estimators such as  kernel and
spline based methods have been
proposed to estimate the regression mean of $Y$ conditional on
$X$. Here the assumption is that
$Y_i=m(X_i)+\epsilon_i$, where for simplicity, we assume $\epsilon_i$
is independent of $X_i$ and has mean zero with density $f_\epsilon(\epsilon)$.
When $X_i$ is unavailable and instead only $W_i$
described above is available, we encounter the problem of nonparametric regression
with measurement error.
It is also well established \citep{fantruong1993} that the same possibly
slow rate of convergence occurs as in the density estimation case.
These convergence rates are also minimax, and again, with a simple
constraint, we can very nearly achieve the non-measurement error
rates.

Working within the Bspline framework, the constraint mentioned above
is that the support of the latent $X$ is compact. This constraint also
arises the in the work of \cite{hall2005}, although they
did not use a Bspline approach.
 Since a Bspline  representation of a function is typically done on a compact support, we use Bspline assisted methods.

\begin{comment}
Our breakthrough came when it struck us that if the functions that we want
to estimate, $f_X(x)$ and
$m(x)$, had been parametric, we would have much simpler problems. If
these problems can be solved while retaining the parameter estimation
rate at the classical root-$n$ as in the error free case,  then
we can approximate $f_X(x)$ and $m(x)$ with Bsplines. Thus,
operationally we would be dealing with parametric models, and there is
a hope that the classical nonparametric convergence rate in the error
free case can also be retained in the measurement error case. A
B-spine representation of a function is typically done on a compact
support, we hence impose this additional assumption. In fact, this
turns out to be the only additional assumption needed
to improve the classical minimax rate in this type of problems.

Subsequently,
we study both density estimation and regression problems in this work
with the help of a spline
representation.
\end{comment}

Using spline representation in measurement  error models is not entirely new, although it was mostly used in  the Bayesian framework \citep{berryetal2002,staudenmayeretal2008,sarkar2014}. While the density estimator is relatively easy to devise, see \cite{staudenmayeretal2008}, regression
estimation turns out to be challenging and it most convenient with our novel further semiparametric treatment. A second challenge in both problems is in establishing the convergence rates of the resulting estimators. The common obstacle in both estimators is the fact that it is a latent function that needs to be approximated with a spline representation, which requires unusual treatment different from the typical handling of the spline approximation. In addition, for the regression problem, we encounter further difficulties because not only the estimator, but also the estimating equation that generates the estimator, do not have
closed forms. We have made novel use of bilinear operators
\citep{conway1990}, which are very different from typical regression
spline asymptotic analysis \citep{smaetal2015}. The detailed proofs
are  in the Appendix  and in an online Supplement.

In the Bspline approximation, the compact set on
which we perform the estimation is built into the procedure at the
very beginning and we benefit from that throughout the
procedure.  Classically, deconvolution is the most widely used method in
nonparametric measurement error problems.
In the typical deconvolution procedures, the Fourier and
inverse Fourier transformation steps do not take advantage
of the compact set knowledge. Instead,  it automatically estimates these
functions on the whole interval, which is much harder to do. How to
modify the deconvolution procedures to make use of the compact support
of $X$ is an interesting future research topic and is well worth
exploring.

In the following, we construct the estimation procedures for both the probability density function and the regression mean function in Section \ref{sec:method}, and summarize the theoretical properties of our estimators in Section \ref{sec:asym}. We provide simulation studies to demonstrate the properties of the new estimators in Section \ref{sec:simulation}, and illustrate the methods in a data example in Section \ref{sec:example}. The paper is concluded with a discussion in Section \ref{sec:discussion}.

\section{Bspline-assisted estimation procedures}\label{sec:method}

\subsection{Probability density function estimation}\label{sec:density}

To set the notation, we use $f_X(x)$ to denote a generic pdf function
of the random variable $X$, and use $f_{X0}(x)$ to denote the true
pdf that generates the data.
We approximate $f_{X0}(x)$ on its support using Bsplines
\citep{masri2005}. For
simplicity, let the support be $[0,1]$. To ensure that
the density function is nonnegative and integrates to 1, we let the
approximation be
\be
f_X(x,\btheta)\equiv\frac{\exp\{\B\trans_r(x)\btheta\}}{
\int_0^1\exp\{\B\trans_r(x)\btheta\}dx}, \label{denest01}
\ee
where $\B_r(x)$ is a vector of  Bspline basis functions, and
$\btheta$ is the Bspline coefficient
  vector. For reasons of identifiability, we fix the
  first component of $\btheta$ at zero, i.e. $\theta_1=0$, and leave
the remaining components $\btheta_L$ free. Here, $a_L$ denotes the subvector of a
generic vector $\a$
without the first
element.
This is different from the estimator of
  \cite{staudenmayeretal2008}. Then
\bse
f_W(w,\btheta)\equiv\frac{\int_0^1\exp\{\B\trans_r(x)\btheta\}f_U(w-x) dx}
{\int_0^1\exp\{\B\trans_r(x)\btheta\}dx}
\ese
is an approximation to the pdf of $W\equiv X+U$, a surrogate of $X$.
We then perform simple  maximum likelihood estimation (MLE), i.e. we maximize
\bse
&&\sumi\log f_W(W_i,\btheta)\\
&=&\sumi\log\int_0^1\exp\{\B\trans_r(x)\btheta\}f_U(W_i-x) dx-
n\log\int_0^1\exp\{\B\trans_r(x)\btheta\}dx
\ese
with respect to $\btheta_L$  to obtain $\wh\btheta_L$, and then reconstruct
$f_X(x,\wh\btheta)$ and use it  as the estimator for $f_{X0}(x)$, i.e.
$\wh f_{X}(x)=f_X(x,\wh\btheta)$. Here $\wh\btheta=(0,\wh\btheta_L\trans)\trans$.

While the estimation procedure for $f_{X0}(x)$ is extremely simple, it
is not as straightforward to establish the large sample properties of
the estimator. In Section \ref{sec:asym}, we will show that
$\wh{f}_{X}(x)$ converges to $f_{X0}(x)$ at a near-nonparametric rate
under mild conditions.

\subsection{Regression function estimation}\label{sec:regression}
Unlike in the density estimation case, the estimation procedure in the
regression model is much more complex. Without considering that the
number of basis functions will increase with sample size,  the key
observation is that as soon as the nonparametric function $m(x)$ is
approximated with the Bspline representation, the regression function
itself without measurement error is a purely parametric model
\citep{wang2009}, hence the idea behind \cite{tsiatisma2004} can be
adapted. Specifically, we treat the Bspline coefficients
  as parameters of interest, treat the unspecified distribution of
  $X$, $f_X(x)$,  as nuisance parameter, and cast the problem as a semiparametric
  estimation problem. We can then construct the efficient score
  function, which relies on $f_{X0}(x)$. A key observation of
  \cite{tsiatisma2004} is that by replacing $f_{X0}(x)$ with an arbitrary
  working model, the consistency of estimation is
  retained. We now describe the estimation procedure in detail.

First, let
$f^*_X(x)$ be a working pdf of $X$. Of course, $f_X^*(x)$ may not be
the same as $f_{X0}(x)$, that is, $f_X^*(x)$ is possibly misspecified.
Following the practice of density estimation in Section
\ref{sec:density}, assume that the true density function $f_{X0}(x)$ has compact support $[0,
  1]$. Therefore we only need to consider $m(x)$ on $[0,1]$. We
  approximate $m(x)$ using the spline representation $\B\trans_r(x)\bb$.
Define
\bse
\S^*_{\bb}(W, Y, \bb) &=& \frac{\partial}{\partial\bb} \log \int_0^1
f_\epsilon\{Y-\B\trans_r(x)\bb\}f_U(W-x) f^*_X(x)d\mu(x)\\
&=& -\frac{\int_0^1 f'_\epsilon\{Y-\B\trans_r(x)\bb\}f_U(W-x) f^*_X(x) \B_r(x) d\mu(x)}{\int_0^1
f_\epsilon\{Y-\B\trans_r(x)\bb\}f_U(W-x) f^*_X(x)d\mu(x)}.
\ese
As the notation suggests, $\S^*_{\bb}(W, Y, \bb)$ is the score
function with respect to $\bb$ calculated from the joint pdf of $(W,
Y)$ under the working model $f_X^*(x)$ and the spline approximation.
Due to the possible misspecification of $f_X^*(x)$, the mean of $\S^*_{\bb}(W, Y,
\bb)$ is not necessarily zero even if the mean function is exactly
$\B\trans_r(x)\bb$. Therefore simply solving $\sumi \S^*_{\bb}(W_i,
Y_i, \bb) = \0$ may generate an inconsistent estimator.
The idea behind our estimator is to find a function $\a(x,\bb)$ so that
\be\label{eq:integral}
&&E\{\S^*_{\bb}(W, Y, \bb)\mid X\}\\
&=&E\left[
\frac{\int_0^1 \a(x,\bb) f_\epsilon\{Y-\B\trans_r(x)\bb\}f_U(W-x) f^*_X(x) d\mu(x)}{\int_0^1
f_\epsilon\{Y-\B\trans_r(x)\bb\}f_U(W-x) f^*_X(x)d\mu(x)}\mid X\right],\n
\ee
and
then solve for $\bb$ using the estimating equation
\be\label{eq:esteq}
&&\sumi \left[\S^*_{\bb}(W_i, Y_i, \bb)\right.\\
&&\left.-
\frac{\int_0^1 \a(x,\bb) f_\epsilon\{Y_i-\B\trans_r(x)\bb\}f_U(W_i-x) f^*_X(x) d\mu(x)}{\int_0^1
f_\epsilon\{Y_i-\B\trans_r(x)\bb\}f_U(W_i-x) f^*_X(x)d\mu(x)}\right]
=\0. \n
\ee
This will guarantee a consistent estimator of $\bb$ if the mean
function is indeed $\B\trans_r(x)\bb$, because our construction
ensures the left hand side of (\ref{eq:esteq}) has mean zero.
The right hand side of (\ref{eq:integral}) is the
conditional expectation of $\a(X,\bb)$ calculated under the Bspline
approximation and the posited model $f_X^*(x)$, hence we alternatively
write it as $E^*\{\a(X,\bb)\mid Y_i, W_i,\bb\}$.

To solve for $\a(x,\bb)$, we discretize the integral equation
(\ref{eq:integral}). In particular, let $f^*_X(x) =
\sum_{j = 1}^L c_j I (x = x_j ) $, where $x_j$'s are points selected
on $[0,1]$, and $c_j\ge0,
\sum_{j = 1}^L c_j = 1$.
Then
\bse
\S^*_{\bb}(W, Y,\bb)
&=& -\frac{\sum_{j=1}^L \B_r(x_j) f'_\epsilon\{Y-\B\trans_r(x_j)\bb\}f_U(W-x_j) c_j }{\sum_{j=1}^L
f_\epsilon\{Y-\B\trans_r(x_j)\bb\}f_U(W-x_j) c_j}.
\ese
Next, to write out the right hand side of (\ref{eq:integral}) upon discretization,
let $\A(\bb)$ be an
$L\times L$ matrix with its $(i, j)$ entry
\bse
A_{ij}(\bb)&=& \int \frac{f_\epsilon\{y-\B\trans_r(x_j)\bb\}f_U(w-x_j)
  c_j}{\sum_{j = 1}^L f_\epsilon\{y-\B\trans_r(x_j)\bb\}f_U(w-x_j) c_j
}\\
&&\hskip 10mm \times f_\epsilon\{y-\B\trans_r(x_i)\bb\}f_U(w-x_i) d\mu(y) d\mu(w).
\ese
Let $\a_i=\a(x_i,\bb)$, $\a=(\a_1, \dots, \a_L)$.
Further, define $\H_j(\bb) = \{H_{1j}(\bb), \ldots, H_{Lj}(\bb)\}$, where
\bse
H_{ij}(\bb)
&=& - \int \frac{
  f'_\epsilon\{y-\B\trans_r(x_j)\bb\}f_U(w-x_j)c_j}{\sum_{j=1}^Lc_j
  f_\epsilon\{y-\B\trans_r(x_j)\bb\}f_U(w-x_j)
  }\\
&&\hskip 10mm \times f_\epsilon\{y-\B\trans_r(x_i)\bb\}f_U(w-x_i)d\mu(y)d\mu(w),
\ese
and let $\b(\bb)$ be a $p\times L$ matrix, with its $i$th column
\bse
\b_i(\bb)
=\sum_{j=1}^L H_{ij}(\bb)\B_r(x_j).
\ese
Then $\b (\bb) = \sum_{j=1}^L
\B_r(x_j)\H_j(\bb)
$.

With this notation, the integral equation (\ref{eq:integral})  is
equivalently written as
$\sum_{j=1}^LA_{ij}\a_j =\sum_{j=1}^L H_{ij}(\bb)\B_r(x_j)$ for $i=1,
\dots, L$, or more concisely,
$\a\A\trans(\bb)=\sum_{j=1}^L
\B_r(x_j)\H_j(\bb)$. Therefore
\bse
\a(\bb)
=\sum_{k=1}^L\B_r(x_k)\H_k(\bb)\{\A^{-1}(\bb)\}\trans,
\ese
hence
\bse
\a_j(\bb)
=\sum_{k=1}^L\B_r(x_k)\H_k(\bb) \{\A^{-1}(\bb)\}\trans\e_j,
\ese
where $\e_j$ is a length $L$ vector with the $j$th component 1 and all
others zero.

Thus, we have obtained $\a(X,\bb)$ on the discrete set $x_1, \dots x_L$ and
can form
\bse
&&E^*\{\a(X,\bb) \mid Y, W, \bb\}\\
&=&\frac{\sum_{j=1}^L \a_j(\bb) f_\epsilon\{Y-\B\trans_r(x_j)\bb\}f_U(W-x_j) c_j}{\sum_{j=1}^L
  f_\epsilon\{Y-\B\trans_r(x_j)\bb\}f_U(W-x_j)
  c_j  }\\
&=&\sum_{k=1}^L\B_r(x_k)\frac{\H_k(\bb) \sum_{j=1}^L \{\A^{-1}(\bb)\}\trans\e_j f_\epsilon\{Y-\B\trans_r(x_j)\bb\}f_U(W-x_j) c_j}{\sum_{j=1}^L
  f_\epsilon\{Y-\B\trans_r(x_j)\bb\}f_U(W-x_j)
  c_j  }\\
&=&\sum_{k=1}^L\B_r(x_k)P_k(W, Y, \bb),
\ese
where
\bse
P_k(W, Y,\bb)=
\frac{\H_k(\bb) \sum_{j=1}^L \{\A^{-1}(\bb)\}\trans\e_j f_\epsilon\{Y-\B\trans_r(x_j)\bb\}f_U(W-x_j) c_j}{\sum_{j=1}^L
  f_\epsilon\{Y-\B\trans_r(x_j)\bb\}f_U(W-x_j)
  c_j  }.
\ese
We then obtain the estimator for $\bb$ by solving  the estimating
equation (\ref{eq:esteq}) with the corresponding $\a(x,\bb)$ plugged in.
In all the functions that are explicitly written to depend on $\bb$, the dependence is
always through $\B_r(\cdot)\trans\bb$.

We show that $\B\trans_r(x)\wh\bb$ converges to $m(x)$ at a nonparametric
rate and we derive its estimation variance in Section \ref{sec:asym}
under mild conditions.
Instead of
  adopting a working model $f_X^*(x)$, we could in fact estimate
  $f_{X0}(x)$ using the method developed in Section \ref{sec:density}
to obtain $\wh f_{X}(x)$, and then use $\wh f_X(x)$ instead of
$f_X^*(x)$ to carry out the subsequent operations. Although
this is the optimal thing to do in theory, we find the practical
improvement of estimation efficiency very small. In addition, unless the
sample size is very large, numerical
instability can arise and sometimes does so. For these reasons, we do
not recommend this approach.

\section{Asymptotic results}\label{sec:asym}

\subsection{Results of probability density function estimation}

We assume the following regularity conditions.

\begin{enumerate}[label=(C\arabic*)]
\item \label{ass:f} The true density function $f_{X0}(x)$ has compact support $[0,
  1]$, is bounded and positive on its support and satisfies
$f_{X0}(x) \in C^q ([0, 1])$, $q\ge1$. The error density $f_U(u)$ is bounded. The conditional
density of $X$ given $W$, $f_{X|W}(x|w)$,  is
bounded. 

\item\label{ass:splineorder}
 The spline order $r \geq  q$.
\item \label{ass:Bknots}
Define the knots $t_{-r + 1}=\dots= t_0 = 0<t_1<\dots<t_N<1=t_{N+1}=\dots=
t_{N+ r}$, where $N$ is the
  number of interior knots and $[0, 1]$ is divided into $N+1 $
  subintervals.  Let $d_{\btheta} = N + r$.
 $N$ satisfies $N \rightarrow
\infty$, $N^{-1} n(\log n )^{-1} \to \infty $ when
 $n\rightarrow \infty$.

\item \label{ass:Bdistance}
 Let $h_p$ be the
  distance between the $(p+1)$st and $p$th
  interior knots and let $h_b = \max\limits_{r\leq p\leq N+r} h_{p}$,
$h_s = \min\limits_{r\leq p\leq N+r} h_{p}$.
  There
  exists a constant $c_{h_b}$, $0 <c_{h_b} < \infty$, such that
\bse
h_b /h_s < c_{h_b}.
\ese
Therefore, $h_b= O(N^{-1})$, $h_s=O(N^{-1})$.

\item\label{ass:dis}
$\btheta_0$ is a spline coefficient with first component $\theta_{01} = 0$ such that $\sup_{x \in [0,
    1]}|\log\{f_X(x, \btheta_0)\}
- \log \{f_{X0}(x)\}| = O_p(h_b^q)$, where $f_X(x, \btheta_0)$ is defined in (\ref{denest01}).
\begin{Rem}$\btheta_0$ exists because by \cite{de1978} there is a
  $\wt{\btheta}$ such
  that
$\sup_{x \in [0,
    1]}|\B_r\trans(x) \wt{\btheta}
- \log \{f_{X0}(x)\}| = O_p(h_b^q)$. Then $\btheta_0 = \wt{\btheta} -
\wt{\theta}_1{\bf 1}$ satisfies the condition, 
  where ${\bf 1}$ is a vector of ones. The argument is facilitated by
  the fact that 
$\B_r\trans(x) {\bf 1}=1$.
\end{Rem}

\item\label{ass:unique}
The expectation
\bse
E\left( \left[\frac{\int_0^1\exp\{\B\trans_{r}(x){\btheta}\}f_U(W_i-x) \B_{rL}(x)
dx }{\int_0^1\exp\{\B\trans_r(x){\btheta}\}f_U(W_i-x) dx }
-\frac{ \int_0^1\exp\{\B\trans_r(x){\btheta}\}\B_{rL}(x) dx
}{\int_0^1\exp\{\B\trans_r(x){\btheta}\}dx}\right]\right)
\ese
is a smooth function of $\btheta_L$ and has unique root in the neighborhood of $\btheta_{0L}$.
\end{enumerate}

Condition \ref{ass:f} contains some basic boundedness and smoothness
conditions on various densities and is quite standard. The only
requirement that appears nonstandard is that $f_{X0}(x)$ has compact
support. This however is practically relevant in many situations as
the range of the values a random variable can be may very well be bounded.
Condition \ref{ass:splineorder} is a standard requirement to ensure that
splines of sufficiently high order are utilized. Condition \ref{ass:Bknots}
requires a suitable amount of spline basis to be used according to the
sample size, and Condition \ref{ass:Bdistance} makes sure that the
spline knots are distributed sufficiently evenly. In summary,
Conditions \ref{ass:splineorder}, \ref{ass:Bknots} and
\ref{ass:Bdistance} are standard requirements and together with
Condition \ref{ass:f}, they ensure
Condition \ref{ass:dis}. We list Condition \ref{ass:dis} instead of
stating it in the proof for convenience. Finally, Condition
\ref{ass:unique} ensures that we are not in the degenerate case where
the expression in \ref{ass:unique} is constantly zero. Under these
conditions, we obtain the following properties of $\wh f_{X}(x)$.

\begin{Pro}\label{pro:consist} Assume Conditions \ref{ass:f} --
  \ref{ass:unique}.
Let $\theta_1 = 0$, and let $\wh{\btheta}_L$ maximize
\bse
&&\sumi\log f_W(W_i,\btheta)\\
&=&\sumi\log\int_0^1\exp\{\B\trans_r(x)\btheta\}f_U(W_i-x) dx-
n\log\int_0^1\exp\{\B\trans_r(x)\btheta\}dx
\ese
with respect to $\btheta_L$,
then $\wh{\btheta}_L - \btheta_{0L} = o_p(1)$.
\end{Pro}

\begin{Pro} \label{pro:theta}
Assume Conditions \ref{ass:f} --
  \ref{ass:unique}.
Let $\theta_1 = 0$, and let $\wh{\btheta}_L$ maximize
\bse
&&\sumi\log f_W(W_i,\btheta) \\
&=&\sumi\log\int_0^1\exp\{\B\trans_r(x)\btheta\}f_U(W_i-x) dx-
n\log\int_0^1\exp\{\B\trans_r(x)\btheta\}dx
\ese
with respect to $\btheta_L$.
Then $\|\wh{\btheta}_L - \btheta_{0L}\|_2 = O_p\{(nh_b)^{-1/2} \}$.
Further, define
\bse
\R_{00}
&=& n^{-1} \sumi\int_0^1\left\{\frac{f_{X0}(x)f_U(W_i-x) }{\int_0^1f_{X0}(x)f_U(W_i-x) dx}- f_{X0}(x)\right\}\B_{rL}(x)dx.
\ese
Then
\bse
&&\wh{\btheta}_L -\btheta_{0L} \\
&=&  \left(E \left[\frac{\int_0^1 f_{X0}(x)f_U(W_i-x) \B_{rL}(x)\B\trans_{rL}(x)
dx }{\int_0^1 f_{X0}(x)f_U(W_i-x) dx}\right.\right.\\
&&\left. \left. - \frac{\left\{\int_0^1f_{X0}(x)f_U(W_i-x) \B_r(x)
dx \right\}^{\otimes2}}{\left\{\int_0^1 f_{X0}(x)f_U(W_i-x) dx \right\}^2}- \int_0^1 f_{X0}(x)\B_{rL}(x)
\B\trans_{rL}(x) dx\right.\right.\\
&&\left.\left.+ \left\{\int_0^1 f_{X0}(x)\B_{rL}(x)
  dx\right\}^{\otimes2}\right]\right)^{-1}\R_{00} \{1 + o_p(1)\}.
\ese
\end{Pro}

Propositions \ref{pro:consist}  and \ref{pro:theta} establishes the
consistency and convergence rate properties of $\wh\btheta_L$ to $\btheta_{L0}$ defined
in Condition \ref{ass:dis}.
The proof of Propositions \ref{pro:consist}  and \ref{pro:theta} are in Supplement
\ref{sec:appconsist} and \ref{sec:apptheta}, respectively. We then
utilize these properties to analyze the bias, variance and
convergence rate of $\wh f_{X}(x)$ in Theorem \ref{th:density}, with its
proof  in Appendix \ref{sec:appdensity}.

\begin{Th}\label{th:density}
Assume Conditions \ref{ass:f} -- \ref{ass:unique}.  Let
$\wh{\btheta}=(0,\wh\btheta_L\trans)\trans$,
$\wh\btheta_L$  be defined in Proposition \ref{pro:theta} and
\bse
\wh f_{X}(x)=\frac{\exp\{\B\trans_r(x)\wh{\btheta}\}}{
\int_0^1\exp\{\B\trans_r(v) \wh{\btheta}\}dv}.
\ese
Then $
\sup_{x \in [0, 1]} |\log\{\wh f_{X}(x)\}
- \log \{f_{X0}(x)\}|
= O_p\{(nh_b)^{-1/2} + h_b^q \}.
$
Specifically,
\bse
{\bias}\{\wh{f}_{X}(x)\}\equiv
E\{\wh{f}_{X}(x)\}-f_{X0}(x)
= O_p(h_b^q)+o_p\{(nh_b)^{-1/2}\}.
\ese
The mean squared error
$\MSE\{\wh f_X(x)\}\equiv\var\{\wh f_X(x)\}+\bias\{\wh
f_X(x)\}^2=O\{h_b^{2q}+(nh_b)^{-1}\}$, and is  minimized at
$N\asymp n^{1/(2q+1)}$ to be $O\{n^{2q/(2q+1)}\}$.
Further,
\bse
\sqrt{nh_b}\left[\wh f_{X}(x)
-f_{X0}(x) - \bias\{\wh{f}(x)\} \right] =1/\sqrt{n}\sumi
Q(W_i,x)+o_p(1).
\ese
If $N n^{-1/(2 q + 1)} \rightarrow \infty$, then
\bse
\sqrt{nh_b}\left\{\wh f_{X}(x)
-f_{X0}(x)  \right\} =1/\sqrt{n}\sumi
Q(W_i,x)+o_p(1).
\ese
Here
\bse
Q(W_i,x)
&=&  \sqrt{h_b}\frac{\partial}{\partial \btheta_L\trans} \frac{\exp\{\B\trans_r(x)\btheta_0\}}
{\int_0^1\exp\{\B\trans_r(x) \btheta_0 \}dx}  \left(E \left[\frac{\int_0^1 f_{X0}(x)f_U(W_i-x) \B_{rL}(x)\B\trans_{rL}(x)
dx }{\int_0^1 f_{X0}(x)f_U(W_i-x) dx}\right.\right.\\
&&- \frac{\left\{\int_0^1f_{X0}(x)f_U(W_i-x) \B_{rL}(x)
dx \right\}^{\otimes2}}{\left\{\int_0^1 f_{X0}(x)f_U(W_i-x) dx \right\}^2}- \int_0^1 f_{X0}(x)\B_{rL}(x)
\B\trans_{rL}(x) dx\\
&&\left.\left.+ \left\{\int_0^1 f_{X0}(x)\B_{rL}(x)
  dx\right\}^{\otimes2}\right]\right)^{-1} \\
&&\times
\int_0^1\left\{\frac{f_{X0}(x)f_U(W_i-x)
  }{\int_0^1f_{X0}(x)f_U(W_i-x) dx}- f_{X0}(x)\right\}
\B_{rL}(x)dx.
\ese
\end{Th}
Theorem \ref{th:density} shows that the Bspline MLE
  density estimator has bias of order $
  O_p(h_b^q)+o_p\{(nh_b)^{-1/2}\} $ and standard error of order $O_p\{(nh_b)^{-1/2}\}$, which is
 the standard
nonparametric density estimation result when no measurement error
occurs, and is a substantial improvement
compared to the minimax convergence rate
without the compact support assumption. In addition, to minimize the
MSE, we can use $h_b\asymp n^{-1/(2q+1)}$, leading to the MSE with order
$n^{2q/(2q+1)}$. On the other hand, to suppress the estimation bias
asymptotically, we need to under-smooth by setting
$h=o\{n^{-1/(2q + 1 )}\}$.

\subsection{Results of regression mean function estimation}

To facilitate the description of the regularity conditions and the
asymptotic results, we introduce some notation. Define
\bse
P(x,W,Y,\bb)&=&\frac{\H(x,\bb) \sum_{j=1}^L \{\A^{-1}(\bb)\}\trans\e_j f_\epsilon\{Y-\B\trans_r(x_j)\bb\}f_U(W-x_j) c_j}{\int_0^1
  f_\epsilon\{Y-\B\trans_r(x)\bb\}f_U(W-x)
 f^*_X(x) d \mu(x)},\\
\H(x,\bb) &=& \{H_{1}(x,\bb), \ldots, H_{L}(x,\bb)\}, \\
H_i(x,\bb)
&=& - \int \frac{
  f'_\epsilon\{y-\B\trans_r(x)\bb\}f_U(w-x)f_X^*(x)}{\int_0^1
  f_\epsilon\{y-\B\trans_r(x)\bb\}f_U(w-x_j) f_X^*(x)d\mu(x)
  }\\
&&\hskip 10mm \times f_\epsilon\{y-\B\trans_r(x_i)\bb\} f_U(w-x_i)d\mu(y)d\mu(w),
\ese
and write
\bse
E^*\{\a(X,\bb) \mid Y, W, \bb\}
&=&\int_0^1 \B_r(x)P(x,W, Y,\bb)d\mu(x),\\
\S^*_{\bb}(W, Y,\bb)
&=& \int_0^1 \frac{\B_r(x) f'_\epsilon\{Y-\B\trans_r(x)\bb\}f_U(W-x) f_X^*(x) }{-\int_0^1
  f_\epsilon\{Y-\B\trans_r(x)\bb\}f_U(W-x) f^*_X(x)}d\mu(x).
\ese

We further define  $
\S^*_{\bb}(W_i, Y_i,m)
$,   $E^*\{\a(X, m) \mid Y_i, W_i, m\}$, $P(x,W, Y,m)$, $P_k(W, Y,m)$
to be the resulting quantities
when we replace all the appearance of
$\B_r(\cdot)\trans\bb$ in $\S^*_{\bb}(W_i, Y_i,\bb)$,
$E^*\{\a(X,\bb) \mid Y_i, W_i, \bb\}$, $P(x,W, Y,\bb)$,
and $P_k(W, Y,m)$ by $m(\cdot)$ respectively.
Here $\a(X, m)$ is a function that satisfies
\bse
E[\S^*_{\bb}(W_i, Y_i,m) - E^*\{\a(X ,m) | Y_i, W_i, m\}|X, m]=\0,
\ese
where the last $m$ is used to emphasize that the
calculation of the outside expectation depends on $m$.
These definitions do not conflict with the notation used
before. In fact, some generalize the previous
notation.
We further define $S_m(Y, W, m)$ to be a linear operator on $L_p$
whose value at $s(\cdot) \in L_p$ is
\bse
&&S_m(Y, W, m)(s) \\
&=& \int_0^1\left[-\frac{
  f'_\epsilon\{Y-m(x)\}f_U(W-x) f^*_X(x)  d\mu(x)}{\int_0^1
f_\epsilon\{Y-m(x)\}f_U(W-x) f^*_X(x)d\mu(x)} - P(x, Y, W, m)\right]
s(x) dx.
\ese
The derivative of $S_m(Y, W, m)$ is a bilinear operator on
$L_p\times L_q$ with $1/p + 1/q = 1, 1 \leq p, q\leq \infty$,
whose value at $s \in L_p, v\in L_q$ is
\bse
\frac{\partial S_m(Y, W, m)}{\partial m} (s, v) =\frac{\partial S_m(Y, W,
m + t v) (s)}{\partial t}\bigg|_{t=0}.
\ese
Further, we define a bilinear operator
\bse
S^2_m(Y, W, m)(s, v) = S_m(Y, W, m)(s) S_m(Y, W, m)(v).
\ese
From the above definition,
\bse
E\{S_m(Y, W, m)\}(s) &=& 0,\\
n^{-1}\sumi S_m(Y_i, W_i, m)(s) - E\{ S_m(Y_i, W_i, m)\}(s) &=& o_p(1) \|s\|_p.
\ese
Also,
\bse
n^{-1}\sumi \frac{\partial S_m(Y_i, W_i, m)}{\partial m} (s, v) -
E\left\{\frac{\partial S_m(Y_i, W_i, m)}{\partial m}\right\} (s, v) &=&
o_p(1)\|s\|_p \|v\|_q,\\
n^{-1}\sumi S^2_m(Y_i, W_i, m)(s, v) - E\{ S^2_m(Y_i, W_i, m)\}(s, v)&=& o_p(1)\|s\|_p \|v\|_q,
\ese
where $E\{ S^2_m(Y_i, W_i, m)\}(s, v) = E\{S_m(Y, W, m)(s) S_m(Y, W,
m)(v)\}$  and\\
$E\{ \partial S_m(Y_i, W_i, m)/\partial m\} (s, v) = E[\{
\partial S_m(Y, W, m)/\partial m\} (s,v)]$.

We now list the regularity conditions.
\begin{enumerate}[label=(D\arabic*)]
\item \label{ass:mf} The true density functions
    $f_{X0}(x),
  f_\epsilon(\epsilon)$ are bounded on their
  supports. In addition, the support of $f_{X0}(x)$ is
 compact.

\item\label{ass:2deriv} $E\left\{{\partial S_m(Y, W, m)}/{\partial
      m}\right\} $ is  a  bounded bilinear operator on  $L_2\times L_2$, $L_1\times
  L_\infty$, and $L_\infty\times L_1$.  Also, $E\{ S^2_m(Y_i, W_i, m)\}$ is
 a  bounded bilinear operator on $L_2 \times L_2$.

\item \label{ass:m} $m(x)$ is bounded on $[0,1]$ and it satisfies
$m(x) \in C^q ([0, 1])$, $q\ge1$.

\item\label{ass:mdis} $\bb_0$ is a $d_\bb$ dimensional spline coefficient such that $\sup_{x \in [0,
    1]}|\B\trans_r(x) \bb_0 - m(x)| = O_p(h_b^q)$. The existence
  of such $\bb_0$ has been shown in \cite{de1978}.

\item\label{ass:munique}
The expectation
\bse
E[\S^*_{\bb}(W_i, Y_i,\bb) - E^*\{\a(X, \bb) | Y_i, W_i, \bb\}]
\ese
has a unique root for $\bb$ in
the neighborhood of $\bb_0$.
Its derivative with respect to $\bb$ is a smooth function of $\bb$ in
the neighborhood of $\bb_0$, with
its singular values bounded and bounded
away from zero.
Denote the unique zero as $\bb^*$.
\end{enumerate}

Conditions \ref{ass:mf} and \ref{ass:2deriv} contain boundedness
requirements on pdfs and operators and are not stringent. Condition
\ref{ass:mf} further requires the compact support of the distribution
of $X$. This is similar to the density estimation case and is crucial.
Conditions \ref{ass:m} and
\ref{ass:mdis} are regarding the mean function $m(x)$ and its spline
approximation which are quite standard. Condition \ref{ass:munique}
is a unique root requirement similar to \ref{ass:unique} and is used
to exclude the pathologic case where the estimating equation is
constantly zero.

In the following,
we establish the consistency of the parameter estimation
in Proposition \ref{pro:betaconsis} and
then further analyze its convergence rate in Proposition
\ref{pro:beta}. The results in these propositions are subsequently
 used to further establish the asymptotic properties of the estimator
 of the regression mean function $m(x)$ in Theorem \ref{th:mean}. The
 proofs of both the propositions and the theorem are  in Supplement
 \ref{sec:appbetaconsis}, \ref{sec:appbeta} and Appendix \ref{sec:appmean} respectively.

\begin{Pro}\label{pro:betaconsis}
Assume Conditions \ref{ass:splineorder} -- \ref{ass:Bdistance},
\ref{ass:mf} -- \ref{ass:munique}.  Let $\wh{\bb}$ satisfy
\bse
 \sumi \left[ \S^*_{\bb}(W_i, Y_i,\wh{\bb}) -E^*\{\a(X,\wh\bb) \mid Y_i, W_i, \wh{\bb}\} \right]=\0.
\ese
Then $\wh{\bb} - \bb_0 = o_p(1)$.
\end{Pro}

\begin{Pro}\label{pro:beta}
Assume Conditions \ref{ass:splineorder} -- \ref{ass:Bdistance},
\ref{ass:m} -- \ref{ass:munique}.  Let $\wh{\bb}$ satisfy
\bse
 \sumi \left[ \S^*_{\bb}(W_i, Y_i,\wh{\bb}) -E^*\{\a(X, \wh{\bb}) \mid Y_i, W_i, \wh{\bb}\}\right] =\0.
\ese
Then $\|\wh{\bb} - \bb_0\|_2 = O_p\{(nh)^{-1/2}\}$. Further,
\bse
\wh{\bb} - \bb_0&=& \left\{E\left(\frac{\partial[\S^*_{\bb}(W_i, Y_i, \bb)
    -E^*\{\a(X, \bb) \mid Y_i, W_i,
\bb\}]}{\partial \bb\trans}\bigg |_{\B_r(\cdot)\trans\bb =
m(\cdot)}\right) \right\}^{-1}\\
&&\times\T_{00}\{1 + o_p(1)\},
\ese
where
\bse
\T_{00} = n^{-1}\sumi \left[ \S^*_{\bb}(W_i, Y_i, m) -E^*\{\a(X, m) \mid Y_i, W_i,
m \}\right].
\ese
\end{Pro}

\begin{Th}
\label{th:mean}
Assume Conditions \ref{ass:splineorder} -- \ref{ass:Bdistance},
\ref{ass:mf} -- \ref{ass:munique}. Let $\wh{m}(x) =
\B\trans_r(x)\wh{\bb}$.
 Then $\sup_{x \in [0, 1]} |\wh{m}(x) - m(x)|
= O_p\{(nh_b)^{-1/2} + h_b^q\}
$. Specifically,  $\bias\{\wh{m}(x)\} = E\{\wh{m}(x)\} - m(x) =O_p(h_b^q)+
o_p\{(nh_b)^{-1/2}\}$.
The mean squared error $\MSE\{\wh m(x)\}\equiv \var\{\wh
m(x)\}+\bias\{\wh m(x)\}^2
=O\{h_b^{2q}+(nh_b)^{-1}\}$, and is minimized at $N\asymp
n^{1/(2q+1)}$
to be $O\{n^{2q/(2q+1)}\}$.
Further,
\bse
\sqrt{n h_b}[\wh{m}(x) - m(x) - \bias\{\wh{m}(x)\}]
= n^{-1/2}\sumi Q(W_i,Y_i,x)+o_p(1).
\ese
If $Nn^{-1/(2q + 1)}\to\infty$, then
\bse
\sqrt{n h_b}\{\wh{m}(x) - m(x)\}
= n^{-1/2}\sumi Q(W_i,Y_i,x)+o_p(1).
\ese
Here
\bse
&&Q(W_i,Y_i,x)\\
&=& \sqrt{h_b} \B\trans_r(x)\left[ -\left\{ E\left(\frac{\partial[\S^*_{\bb}(W_i, Y_i, \bb)
    -E^*\{\a(X, \bb) \mid Y_i, W_i,
\bb\}]}{\partial \bb\trans}\bigg |_{\B_r(\cdot)\trans\bb =
m(\cdot)}\right)\right\}^{-1}\right.\\
&&\left.\times  \S^*_{\bb}(W_i, Y_i, m) -E^*\{\a(X, m) \mid Y_i, W_i,
m \} \right].
\ese
\end{Th}

Theorem \ref{th:mean} shows that the Bspline
  regression mean function estimator has bias of order
  $O_p(h_b^q) +o_p\{(nh_b)^{-1/2}\}$ and standard error of order $O_p\{(nh_b)^{-1/2}\}$, which is
 the  standard
nonparametric regression result without measurement error,
 and is much better than the minimax convergence rate without assuming
 $X$ to be compactly supported. Further, to minimize the
MSE, we let $h_q\asymp n^{-1/(2q+1)}$, leading to the MSE of order
$n^{2q/(2q+1)}$. On the other hand, to suppress the estimation bias
asymptotically, we need to under-smooth by setting
$h=o\{n^{1/(2q + 1)}\}$.

\section{Simulation studies}\label{sec:simulation}

\subsection{Performance of the Bspline-assisted density estimator}

\begin{figure}[!h]
\centering
\includegraphics[scale = 0.3]{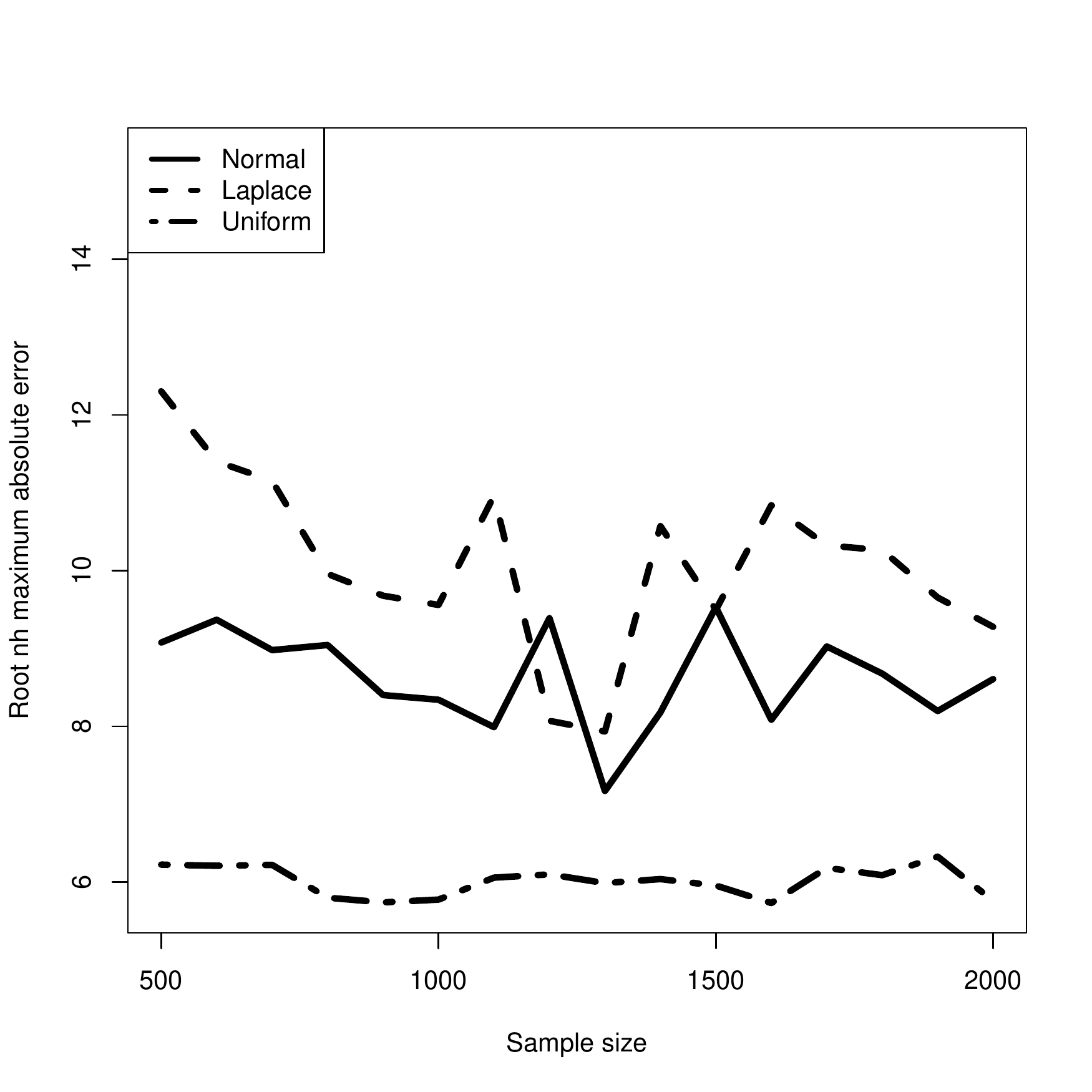}
\includegraphics[scale = 0.3]{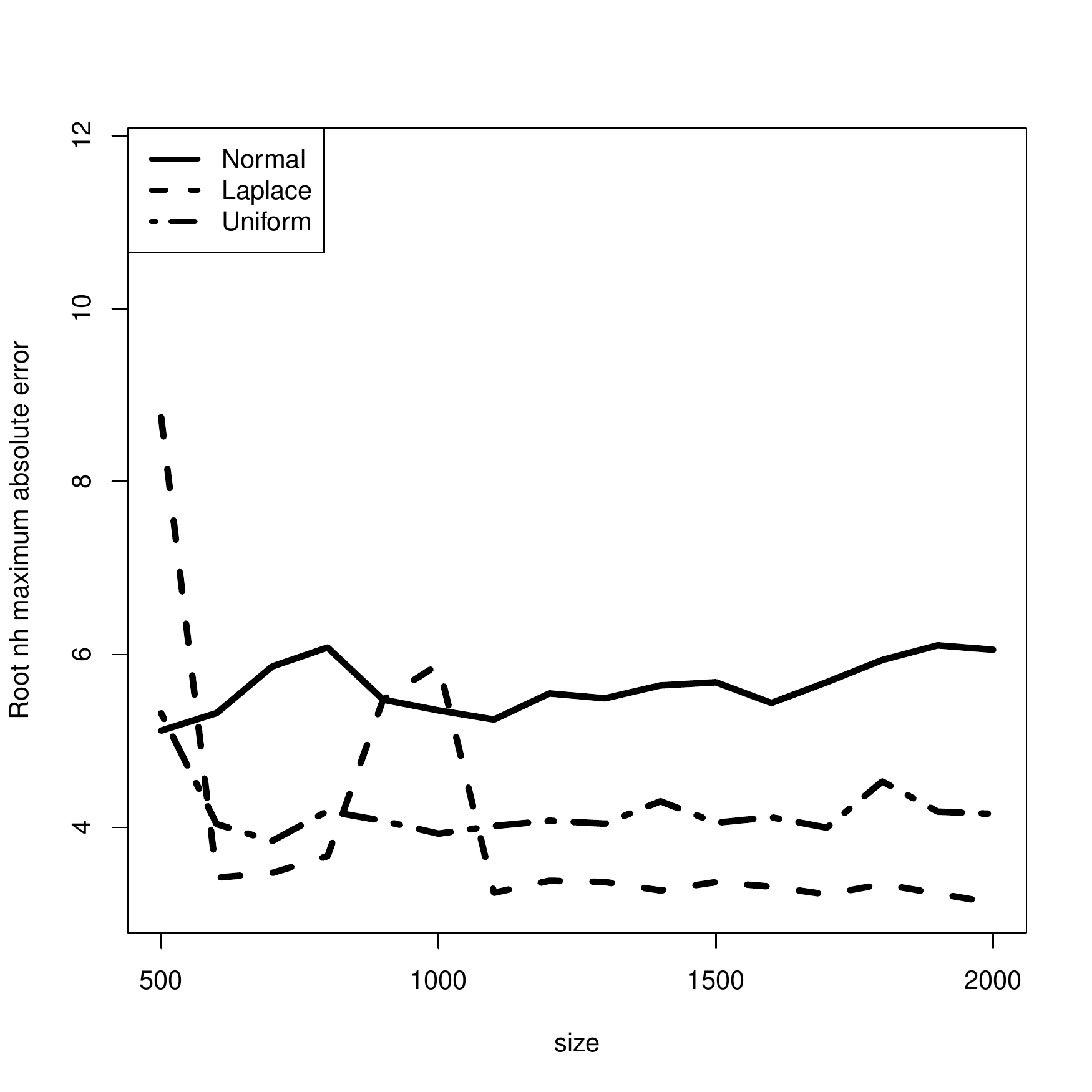}
\caption{Performance of the Bspline MLE pdf estimation (left) and
Bspline semiparametric mean estimation (right). Results based
  on 200 simulations. The solid lines are for the Normal measurement
  error case, the dashed lines are for Laplace error and the
  dot-dashed lines are for the Uniform error case. 
    Change the dotted line for the Laplace to a dashed line, it should
    be Normal, Laplace and Uniform, and you should make the lines
    thicker. Also, add the text "The solid lines are for the Normal
    measurement error case, the dashed lines are for Laplace error and
    the dot-dashed lines are for the Uniform error case."}
\label{fig:simurate}
\end{figure}

We conducted two simulation studies to illustrate the finite sample
performance of the proposed pdf and regression mean estimators.
To evaluate the Bspline MLE method for estimating the density functions,
we generated data sets of sample sizes from $n
= 500$ to $n=2000$. We used cubic Bsplines
with the number of knots equal to the smallest integer larger than $1.3n^{1/5}$.
In all our simulations, 
$X$ is generated from a beta distribution with both shape parameters
equal to 4. We then generated the measurement error $U$ from three models:
\begin{itemize}
\item[]I(a): a normal distribution with mean
0, variance $0.25$, denoted by $N(0, 0.25)$,
\item[]I(b):  a Laplace
distribution with mean 0 and scale $0.5/\sqrt{2}$, denoted by ${\rm Lap}(0,
0.5/\sqrt{2})$,
 \item[]I(c): a uniform distribution on $[-\sqrt{3/4}, \sqrt{3/4}]$, denoted by\\
${\rm Unif} (-\sqrt{3/4}, \sqrt{3/4})$.
\end{itemize}
In the left panel of Figure \ref{fig:simurate}, we plotted the
averaged  root-$(nh_b)$
maximum absolute error (MAE), calculated as $\sqrt{n
h_b} \sup_x|\wh f_{X}(x) - f_{X0}(x) |$, versus the sample sizes $n$.
Following Theorem \ref{th:density}, the root-$(nh_b)$ MAE has a constant
order.

 This translates to the curves in the plots that stabilize in a
small range, which is what we observe, especially when the sample size
grows to larger than 700.

\begin{figure}[!h]
\centering
\caption{Comparison of pdf (upper) and mean (lower) estimators based on the Bspline
  MLE or Bspline semiparametric estimator
  (solid) and the
  deconvolution (dashed) method, when measurement errors are Normal (left),
 Laplace (middle) and Uniform (right) respectively.
Average maximum
 absolute error $\sup_x|\wh f_X(x)-f_X(x)|$
(upper) or $\sup_x|\wh m_X(x)-m_X(x)|$ (lower)
 is computed
based on 200 simulations at sample sizes from 500 to 2000.}
\label{fig:compare}
\includegraphics[scale = 0.23]{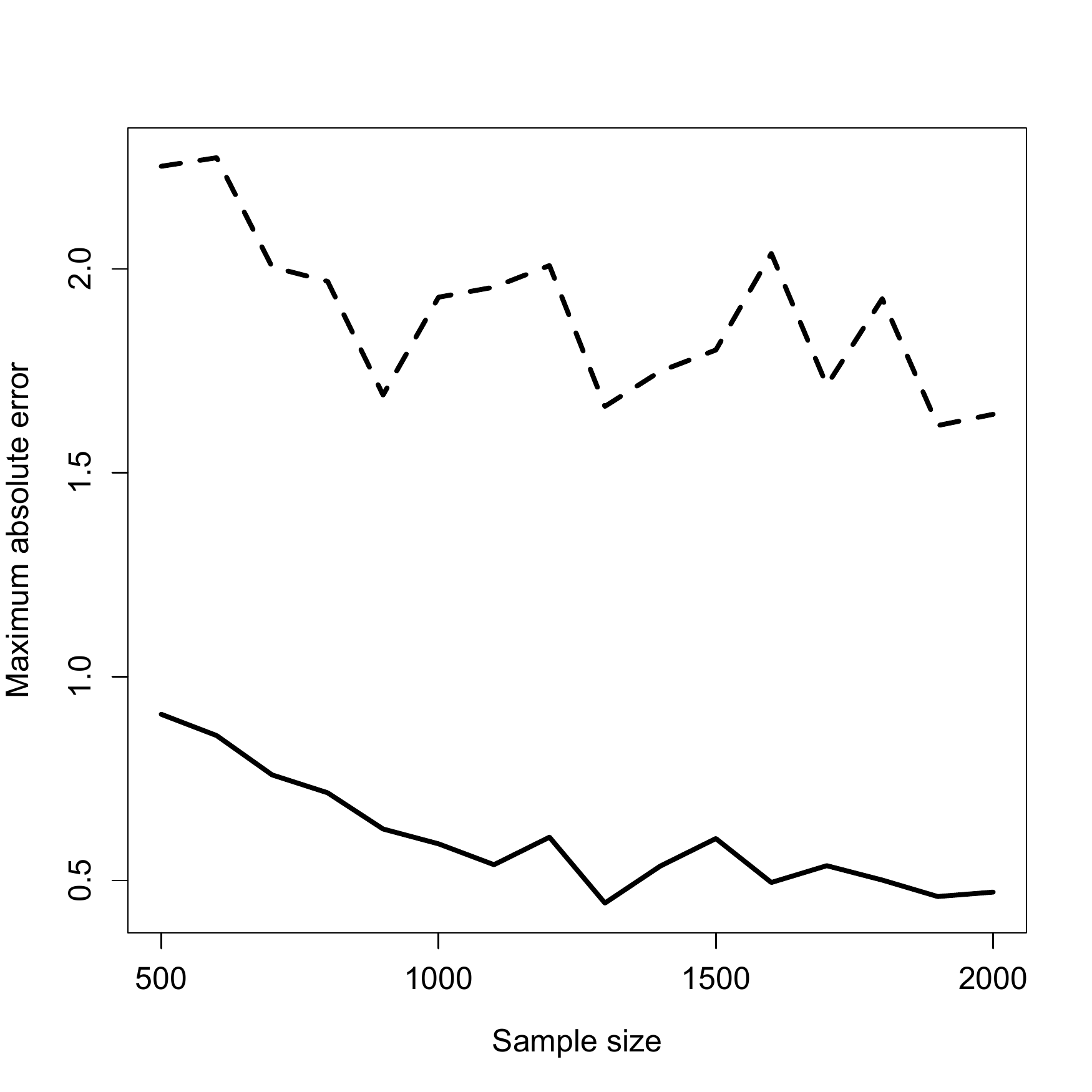}
\includegraphics[scale = 0.23]{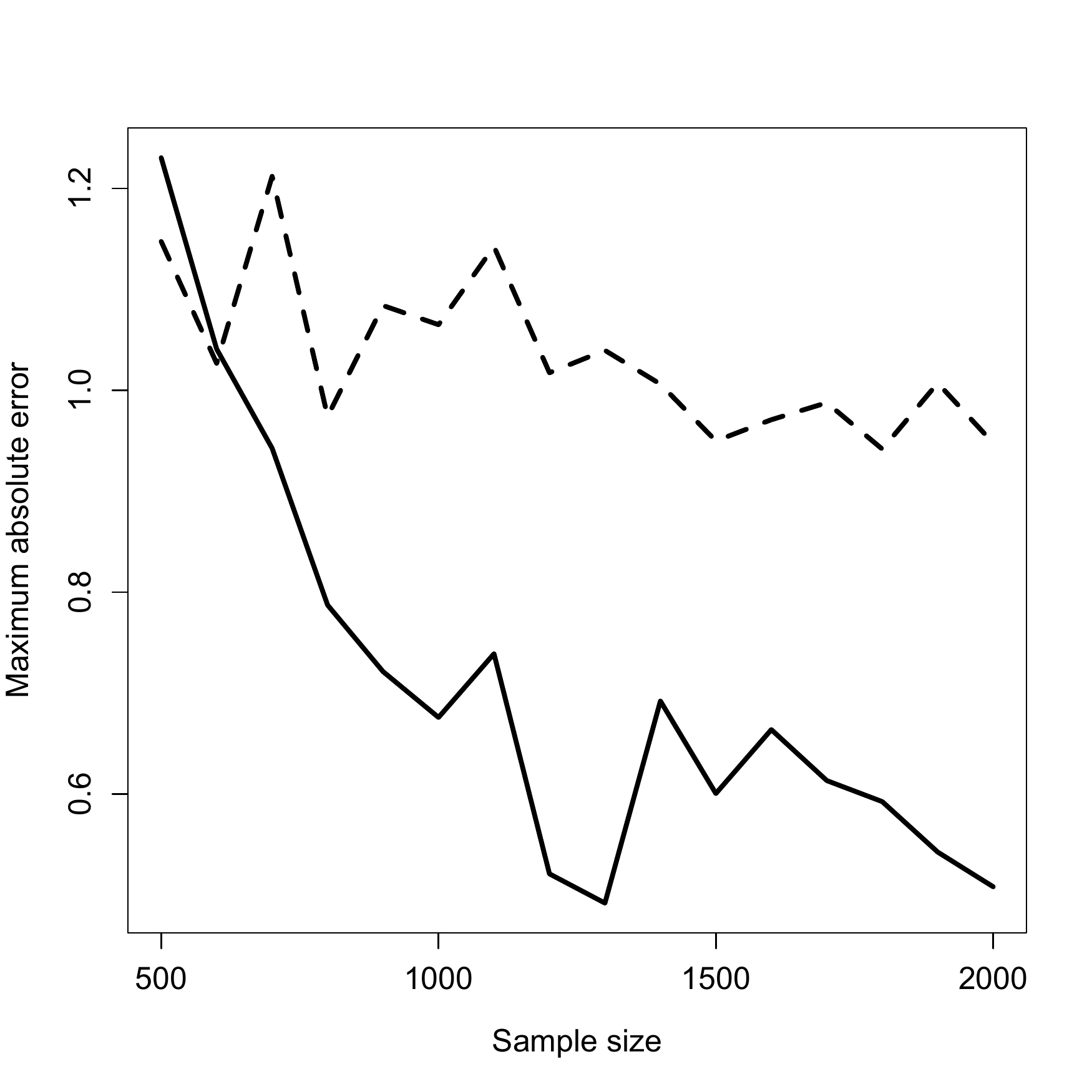}
\includegraphics[scale = 0.23]{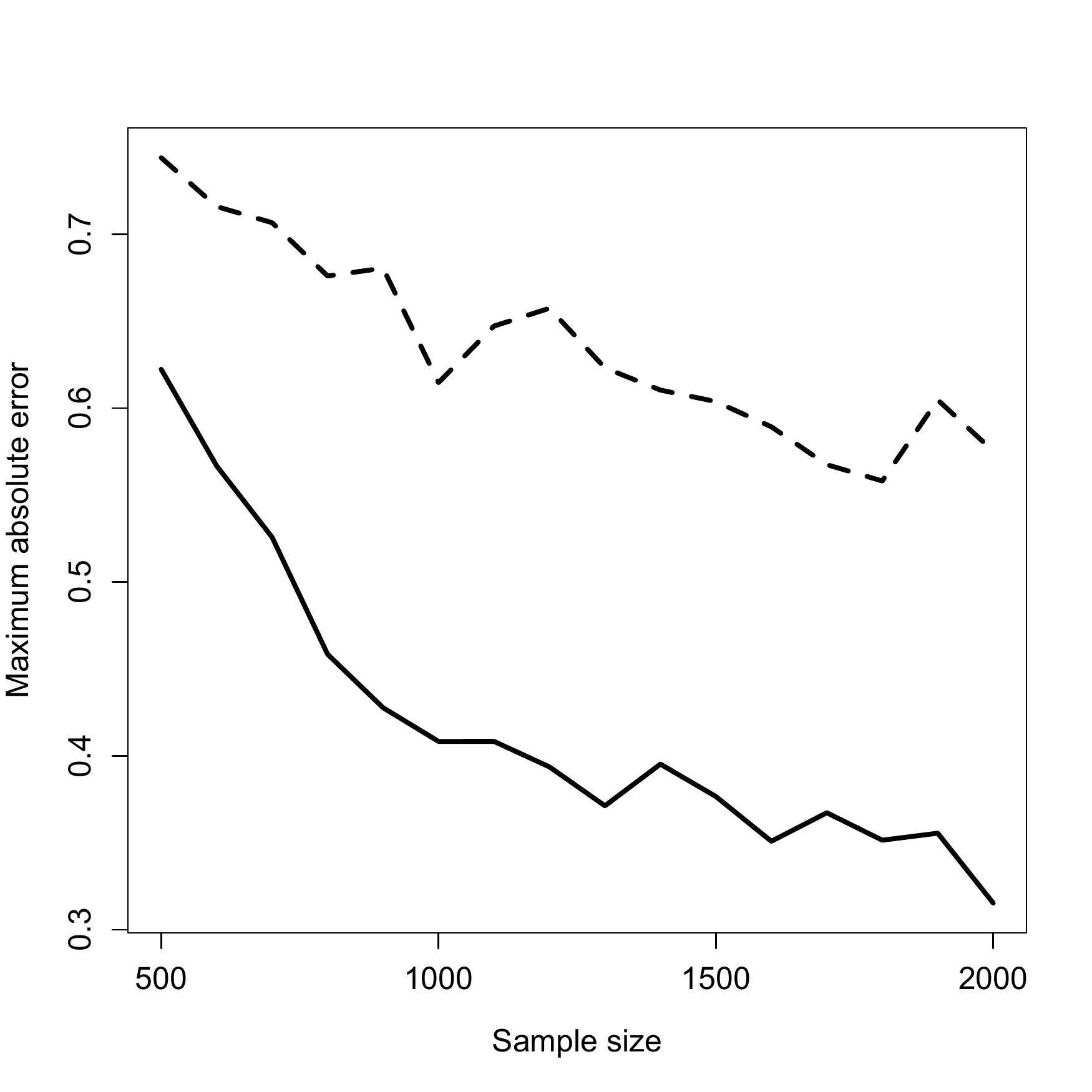}\\
\includegraphics[scale = 0.23]{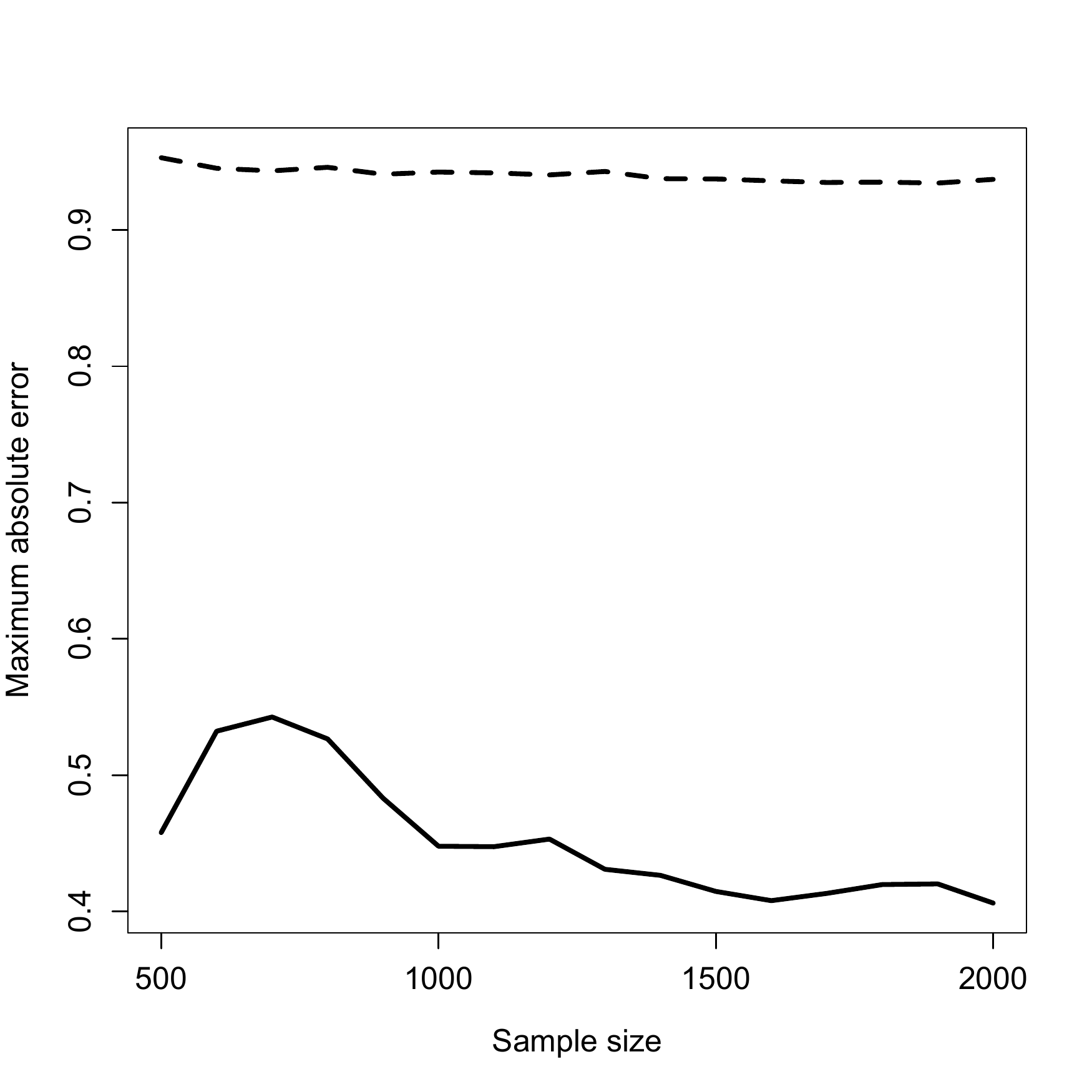}
\includegraphics[scale = 0.23]{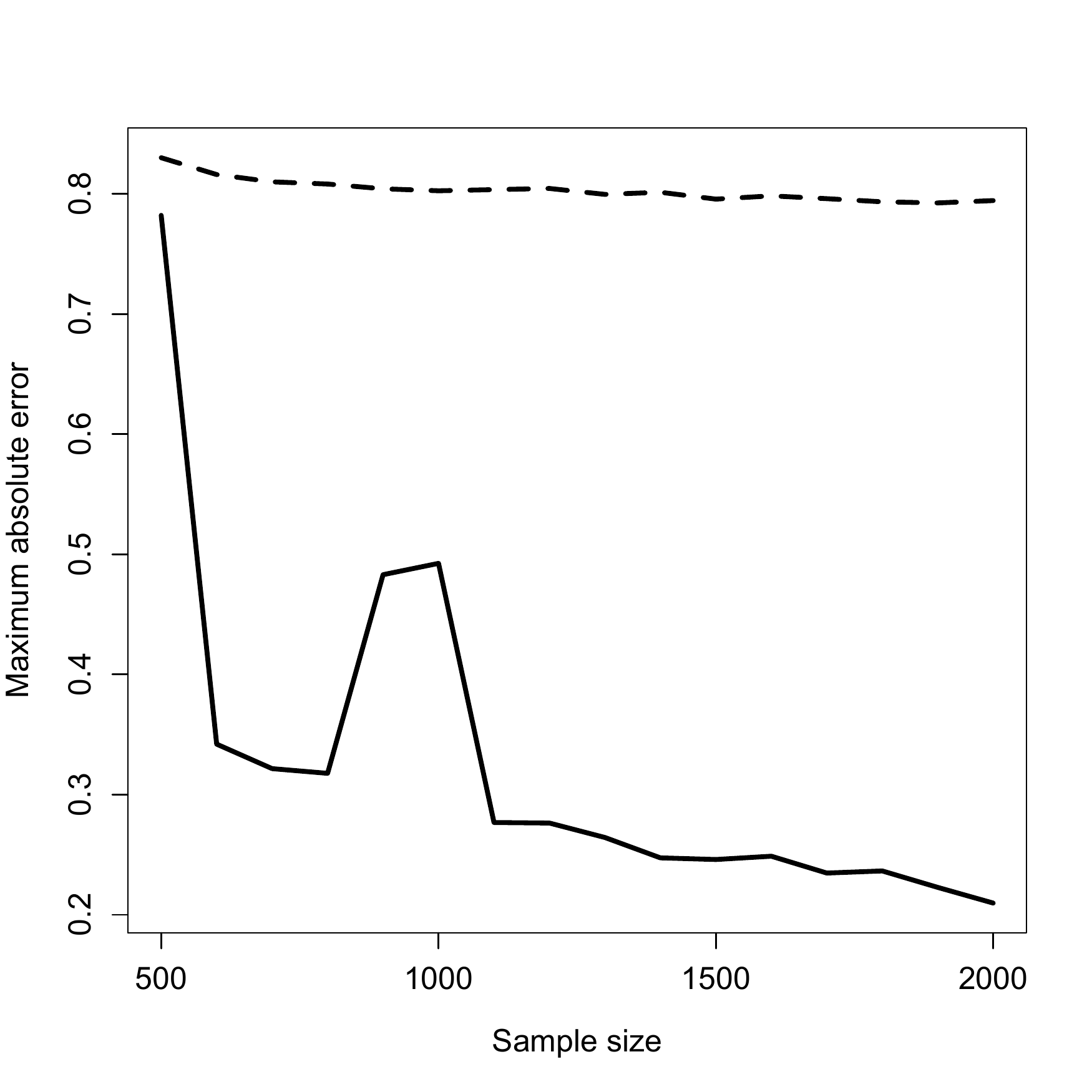}
\includegraphics[scale = 0.23]{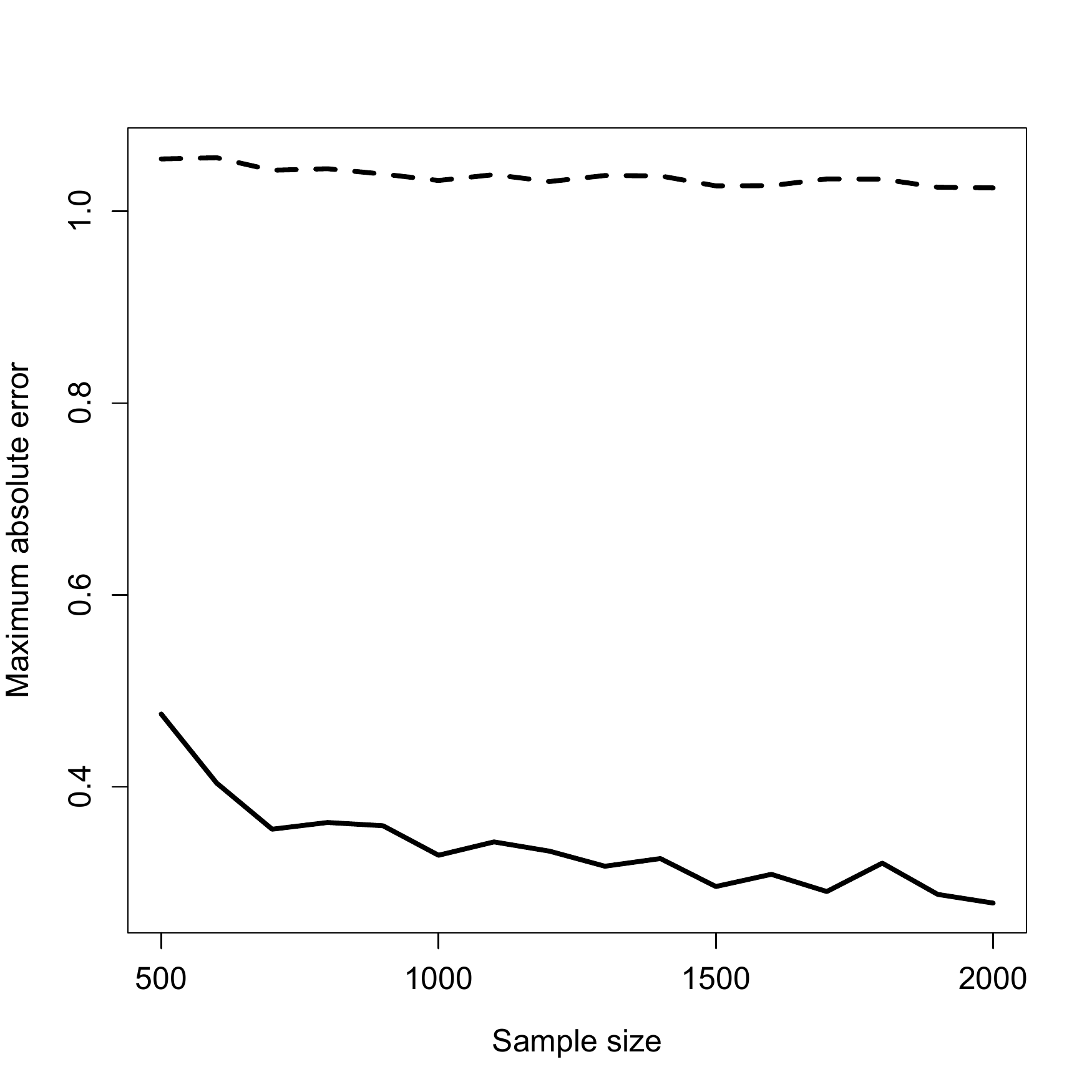}
\end{figure}

We also compared the Bspline MLE method with the widely used
deconvolution method  \citep{stefanskicarroll1990} for
density estimation that does not assume $X$
is compactly supported. In the upper row of
Figure \ref{fig:compare}, we plotted the average values of
$\sup_x|\wh f_{X}(x) - f_{X0}(x)|$  based on 200 simulations for both
methods at different sample sizes. We adopted the two-stage plug-in
bandwidth selection method proposed in \cite{delaigle2002}
in implementing the deconvolution method. Unsurprisingly,
results in the upper row of Figure \ref{fig:compare} indicate that the
Bspline MLE method outperforms the deconvolution method with rather
significant gain in this case. We suspect that this is because the
measurement errors are quite large here which caused difficulties for
both methods, but especially for the deconvolution method.

\begin{figure}[!h]
\caption{Bspline MLE density estimation (left) and deconvolution
  estimation (right) from 200  simulations: The solid
  lines represent the true functions and the
  dashed lines represent the estimated functions and their 90\%
  confidence bands. The first row to third row are the results for
  Model II (a)--(c), respectively. Sample size is
  500. The figures on the right should have the
    y-axis label to be $m$, not $f$.
}
\label{fig:den1}
\begin{center}
\includegraphics[scale = 0.35]{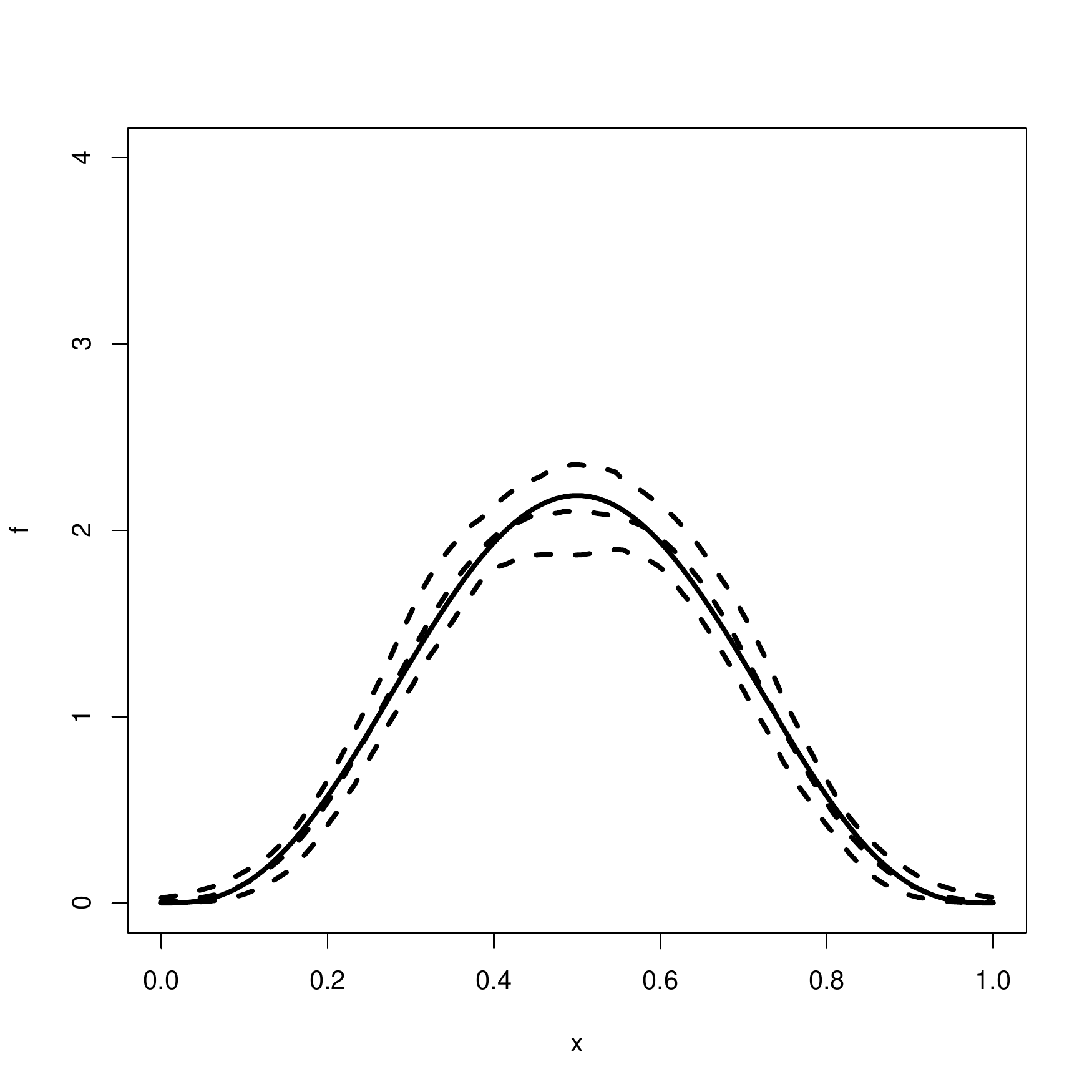}
\includegraphics[scale = 0.35]{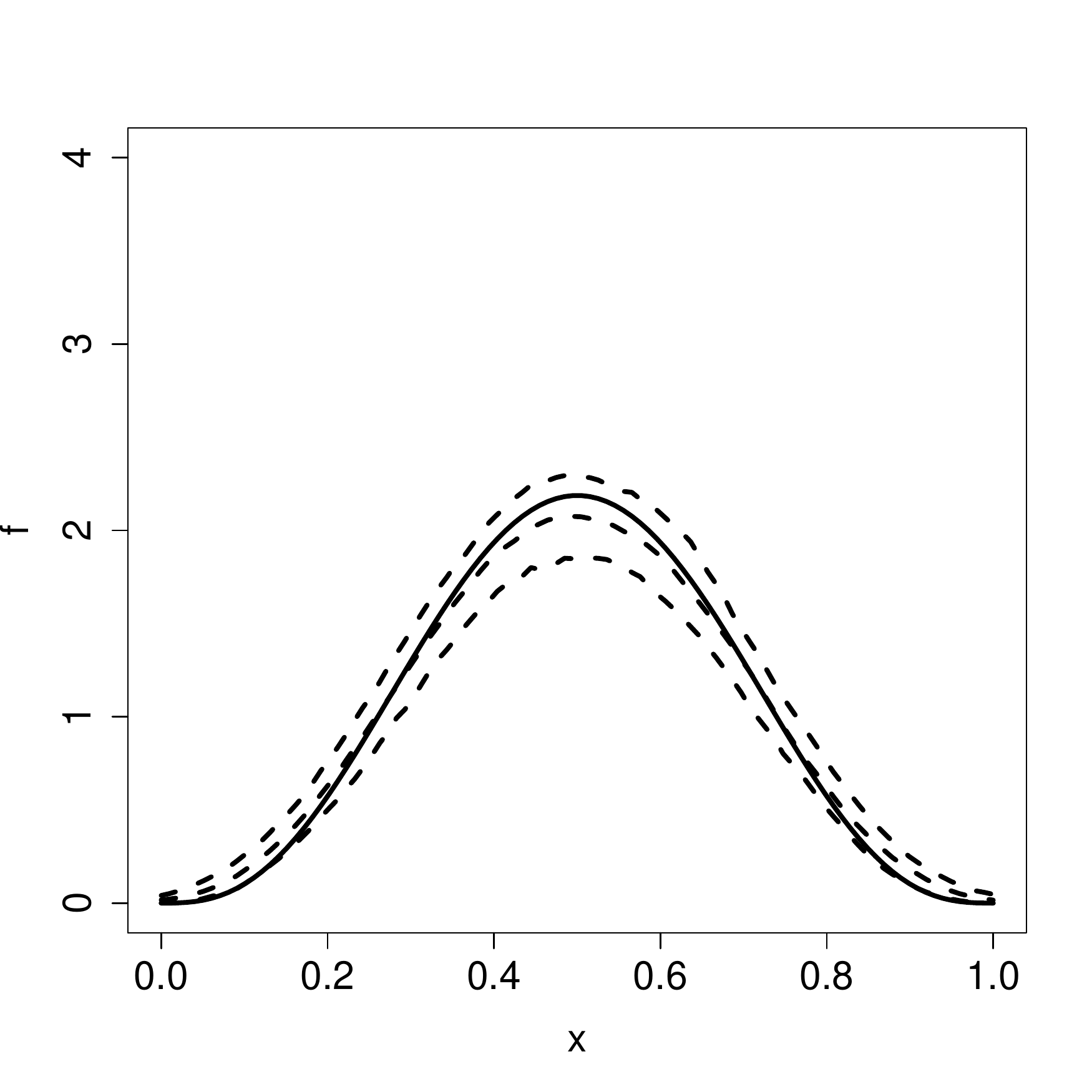}\\
\includegraphics[scale = 0.35]{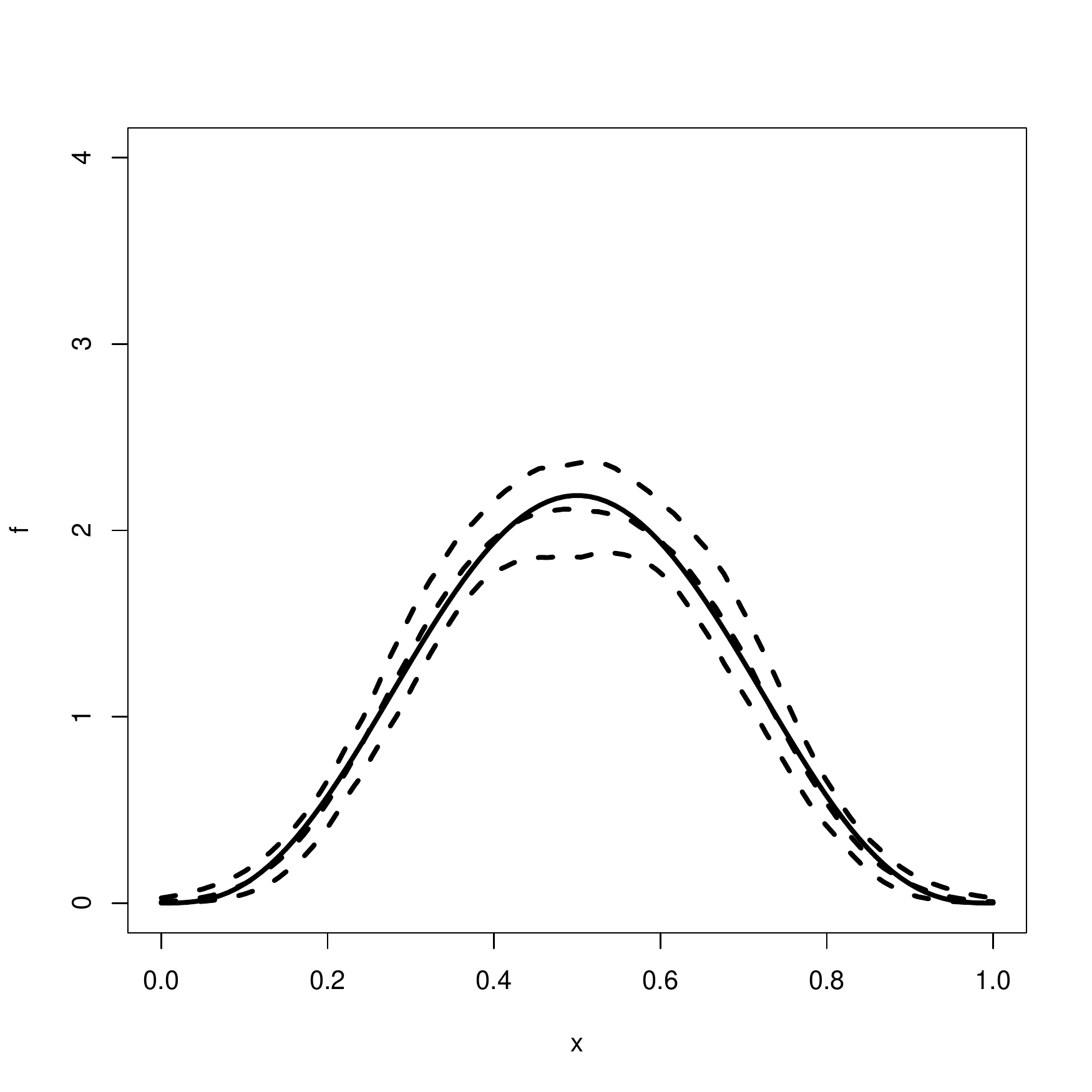}
\includegraphics[scale = 0.35]{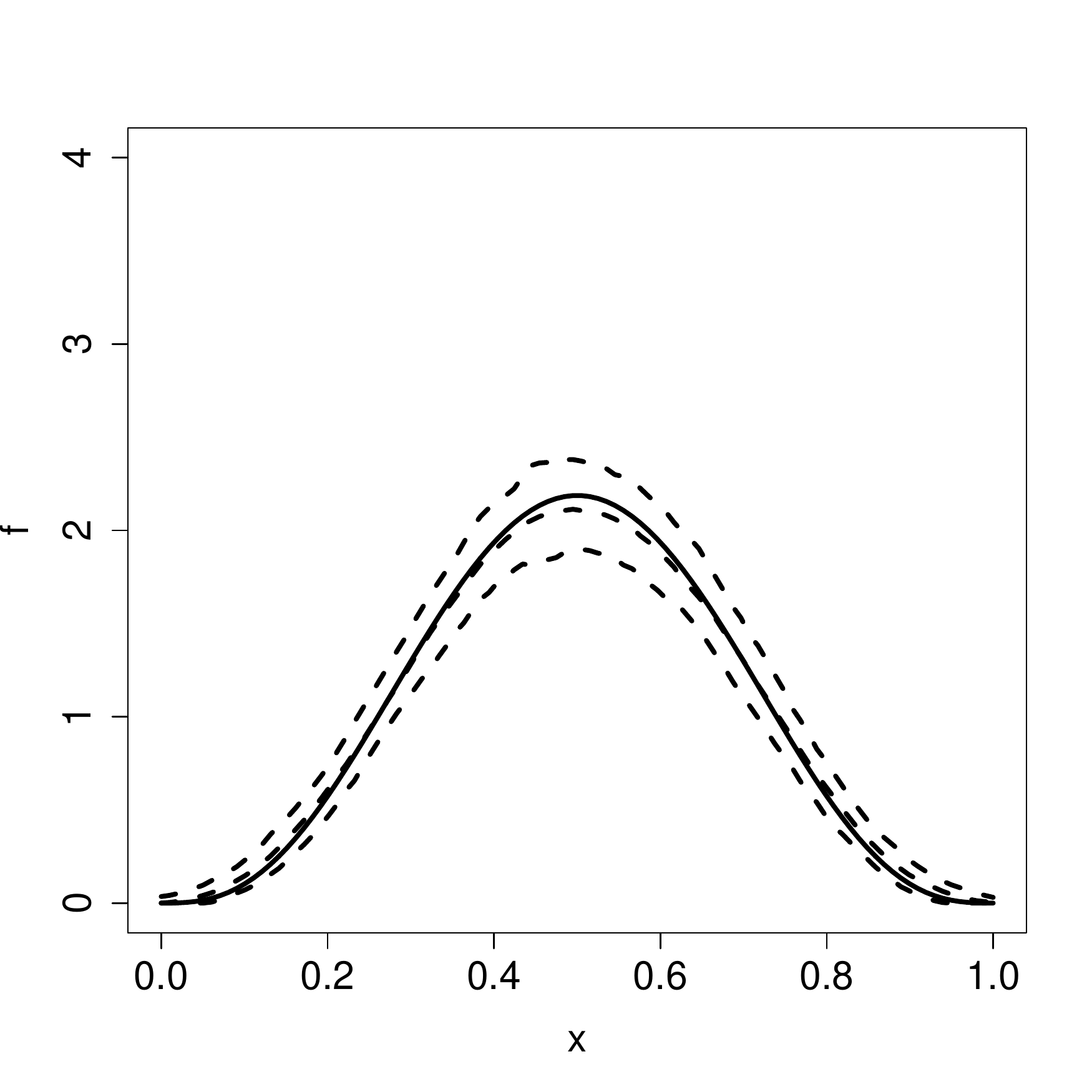}\\
\includegraphics[scale = 0.35]{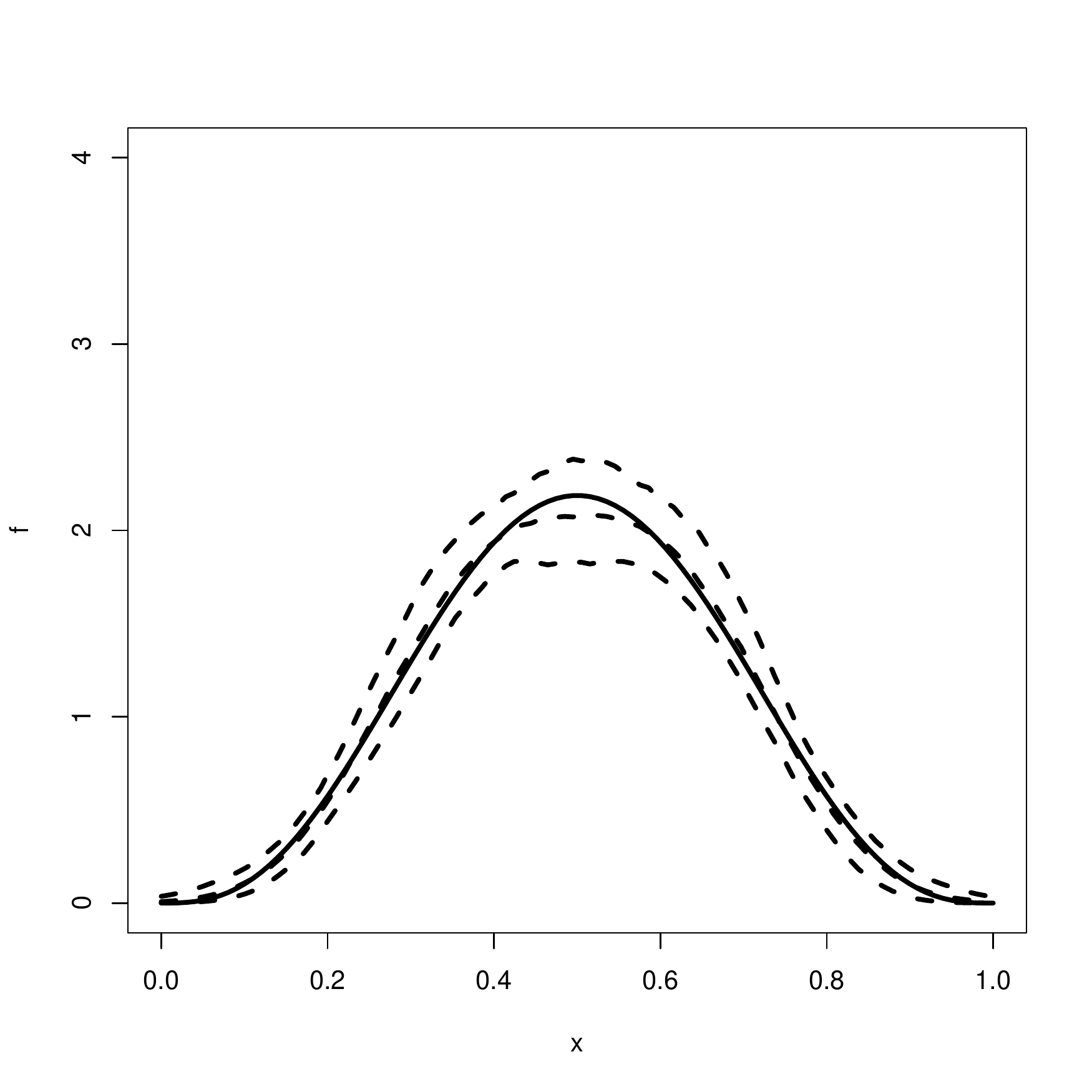}
\includegraphics[scale = 0.35]{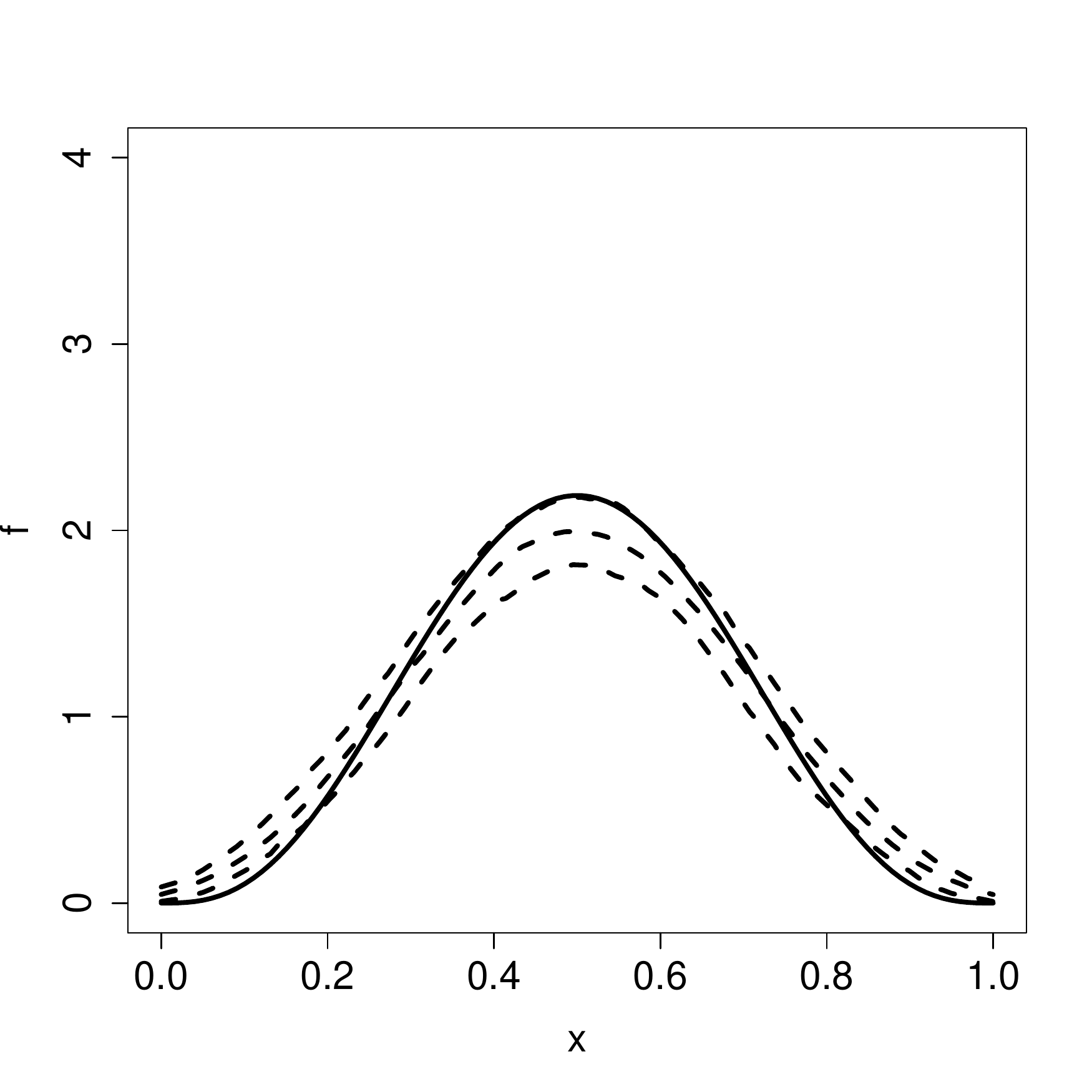}
\end{center}
\end{figure}

To further examine the performance of individual estimated pdf
curves from both methods, we also plotted the estimated mean density curves
for sample sizes 500, 1000
and 2000. Because the deconvolution method performs poorly when the
measurement errors are large (see upper row of Figure \ref{fig:compare}), we
reduced the error variances, and
generated the three error distributions from models
\begin{itemize}
\item[] II(a): $N(0, 0.0025)$,
II(b): ${\rm Lap}(0,
0.05/\sqrt{2})$,
II(c): ${\rm Unif} (-0.125, 0.125)$.
\end{itemize}
Under the reduced error variability, we performed the estimation and
plotted the resulting density estimates and their 90\%
confidence bands.
Figure \ref{fig:den1}  contains the results of the Bspline MLE and
the deconvolution estimator using the
two-stage plug-in bandwidth \citep{delaigle2002}, at the sample size
500. Although not as
dramatic as the large error case,
the Bspline MLE is still closer to the true pdf with narrower confidence
band, hence is more precise than the deconvolution
method.

We also provide similar comparisons at sample sizes
1000 and 2000 in Figures \ref{fig:den2} and \ref{fig:den3} in the
Supplement.
For a quantitative comparison, we also computed the
MAE between the true pdf curve and the estimated
pdf curve in the upper part of Table \ref{tab:compare}. These results
show that the Bspline MLE performs consistently better than
deconvolution which does not make the compact support assumption.
The performance of the two methods largely follow the same pattern
when sample sizes are smaller, although the improvement of the Bspline
method over the deconvolution method is not as dramatic, naturally. We provide
the results for $n=200$ in Figure \ref{fig:den0} in the Supplement.

\begin{comment}
To further examine the performance of individual estimated pdf
curves from both methods, we also plotted the estimated density curves
for sample sizes 200 500, 1000
and 2000. Because the deconvolution method performs poorly when the
measurement errors are large (see upper row of Figure \ref{fig:compare}), we
reduced the error variances, and
generated the three error distributions from models
\begin{itemize}
\item[] II(a): $N(0, 0.0025)$,
II(b): ${\rm Lap}(0,
0.05/\sqrt{2})$,
II(c): ${\rm Unif} (-0.125, 0.125)$.
\end{itemize}
Under the reduced error variability, we performed the estimation and
plotted the resulting density estimates and their 90\%
confidence bands.
Figure \ref{fig:estden} contains the  results of the Bspline MLE and
Figure \ref{fig:cpden} of
 deconvolution estimator using the
two-stage plug-in bandwidth \citep{delaigle2002}. Comparing
Figures \ref{fig:estden} and  \ref{fig:cpden}, although not as
dramatic as the large error case,
the Bspline MLE
 is still closer to the true pdf with narrower confidence
band, hence is more precise than the deconvolution
method. For easier comparison, we also computed the
MAE between the true pdf curve and the estimated
pdf curve in the upper part of Table \ref{tab:compare}. These results
show that the Bspline MLE performs consistently better than
deconvolution.
\end{comment}

\subsection{Performance of the Bspline-assisted semiparametric estimator}
We next evaluated the finite sample performance of the
Bspline semiparametric mean regression method described in Section
\ref{sec:regression}. We used sample sizes  $n$ from
500 to 2000, used cubic Bsplines, with the number of knots
 equal to the smallest integer larger than $1.3n^{1/5}$.
In this case we generated $X$ from a beta distribution with both
shape parameters equal to 2.
The true regression mean function is
$m(x) = \sin(2 \pi x)$
and we generated the regression model errors $\epsilon$
from a normal distribution with mean zero and variance 0.25.
We generated the measurement errors $U$ from the three different
distributions described in Models I(a)--I(c).

 In the right panel of Figure
\ref{fig:simurate},  we plotted the averaged root-$(nh_b)$ MAE
calculated via $\sqrt{n
h_b}\sup_x|\wh{m}(x) - m(x)|$ as a function of the sample size
$n$. Similar to the density estimation experiment, the curves stabilize
as sample size increases, and is largely flat after $n=1000$,
indicating that $\sup_x|\wh{m}(x) - m(x)|$ has order $(nh_b)^{-1/2}$.

We further compared the Bspline semiparametric method with the deconvolution
method \citep{fantruong1993} in the nonparametric mean
regression model, where the deconvolution method does not make the
compact support assumption on $X$.
 In the lower row of
Figure \ref{fig:compare}, we plotted the average $\sup_x|\wh{m}(x) - m(x)
|$  over 200 simulations for both methods. Again, in this case, with
moderate to significant amount of noise, the
Bspline semiparametric method greatly outperforms the deconvolution
method with smaller average error.

\begin{figure}[!h]
\caption{Bspline semiparametric regression estimation (left) and
deconvolution estimation (right)
from 200
  simulations: The solid lines represent the true functions and the
  dashed lines represent the estimated functions and their 90\%
  confidence bands. The rows are the results for
  Models II (a)--(c), respectively. Sample size is
  500. The figures on the right should have the
    y-axis label to be $f$, not $m$.}
\label{fig:mean1}
\begin{center}
\includegraphics[scale = 0.35]{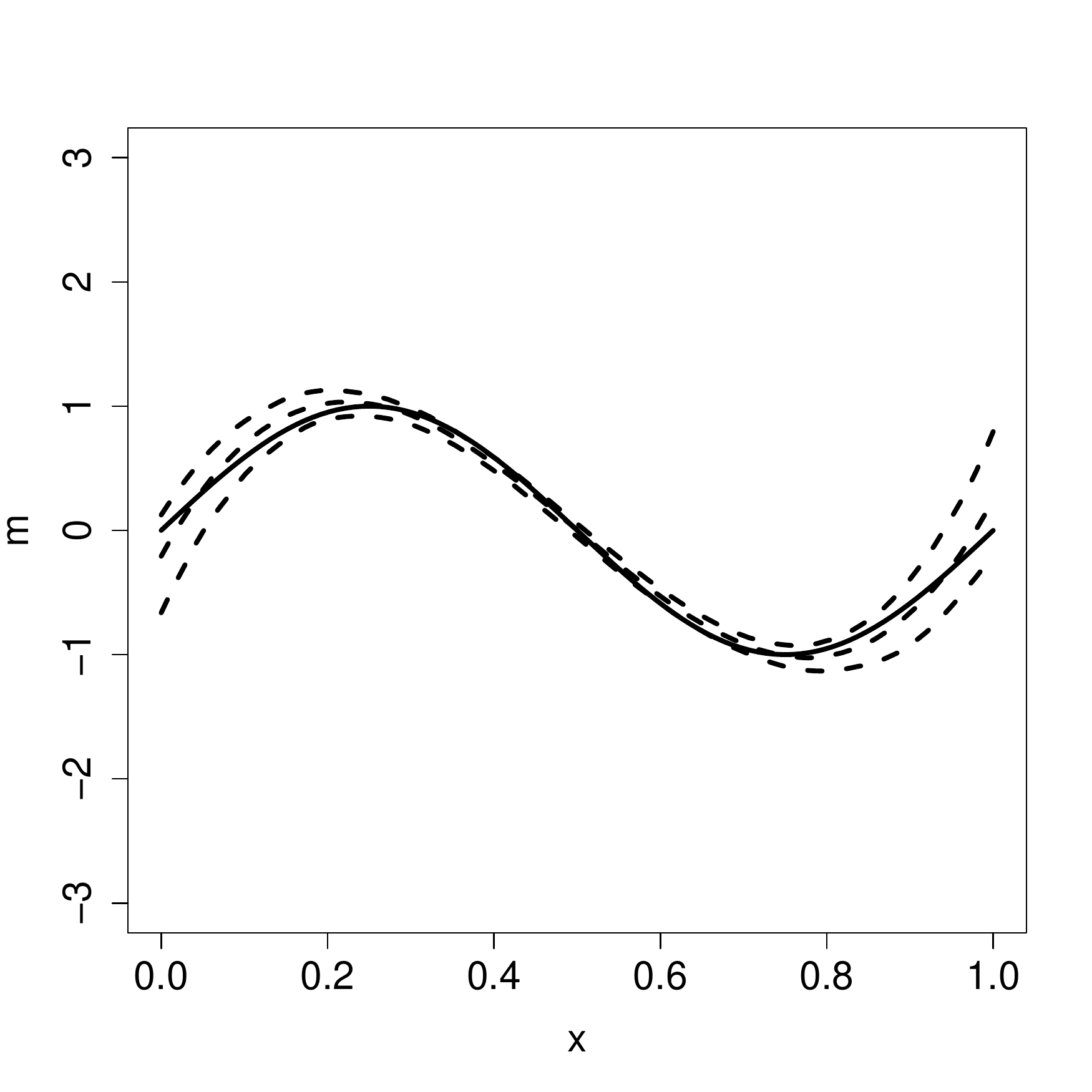}
\includegraphics[scale = 0.35]{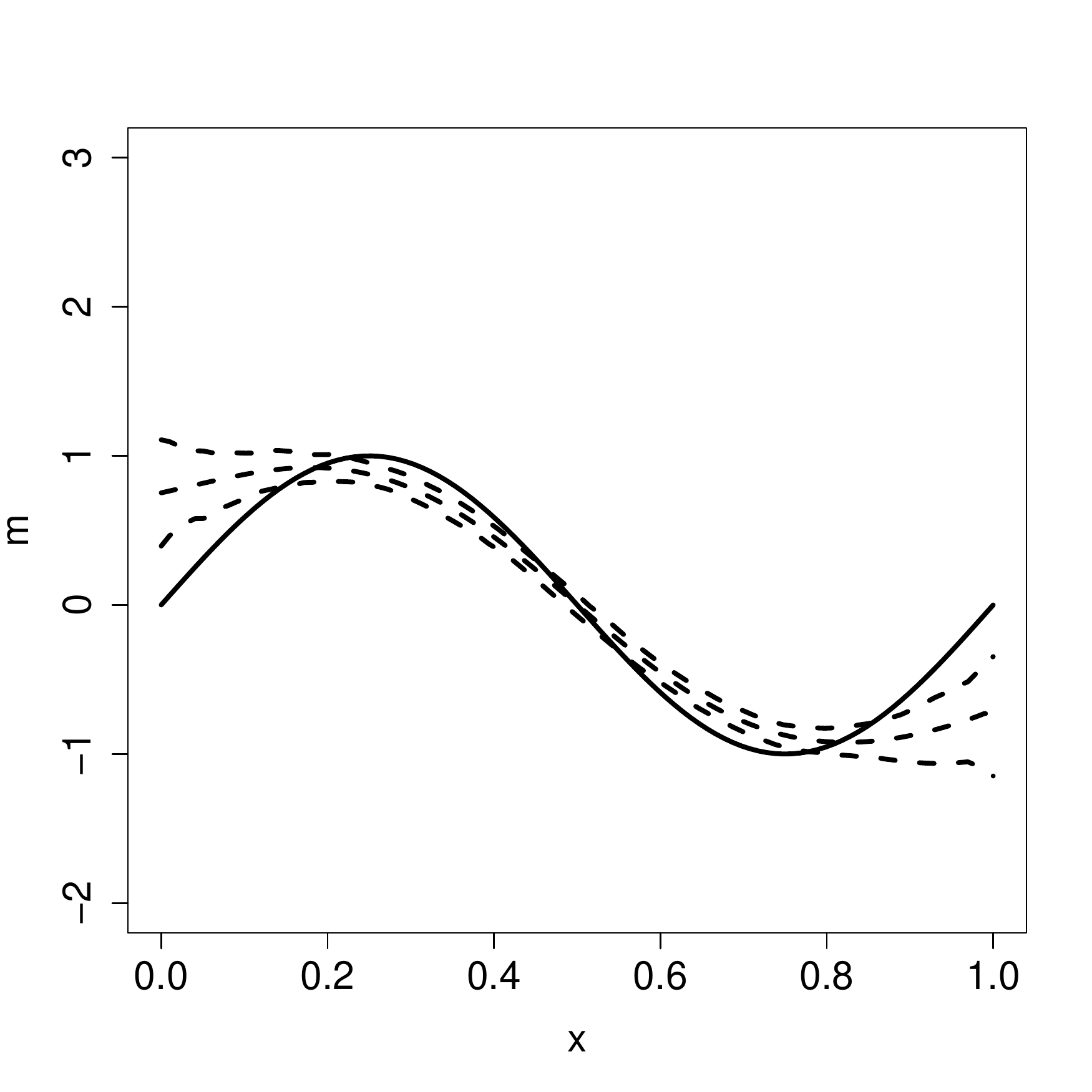}\\
\includegraphics[scale = 0.35]{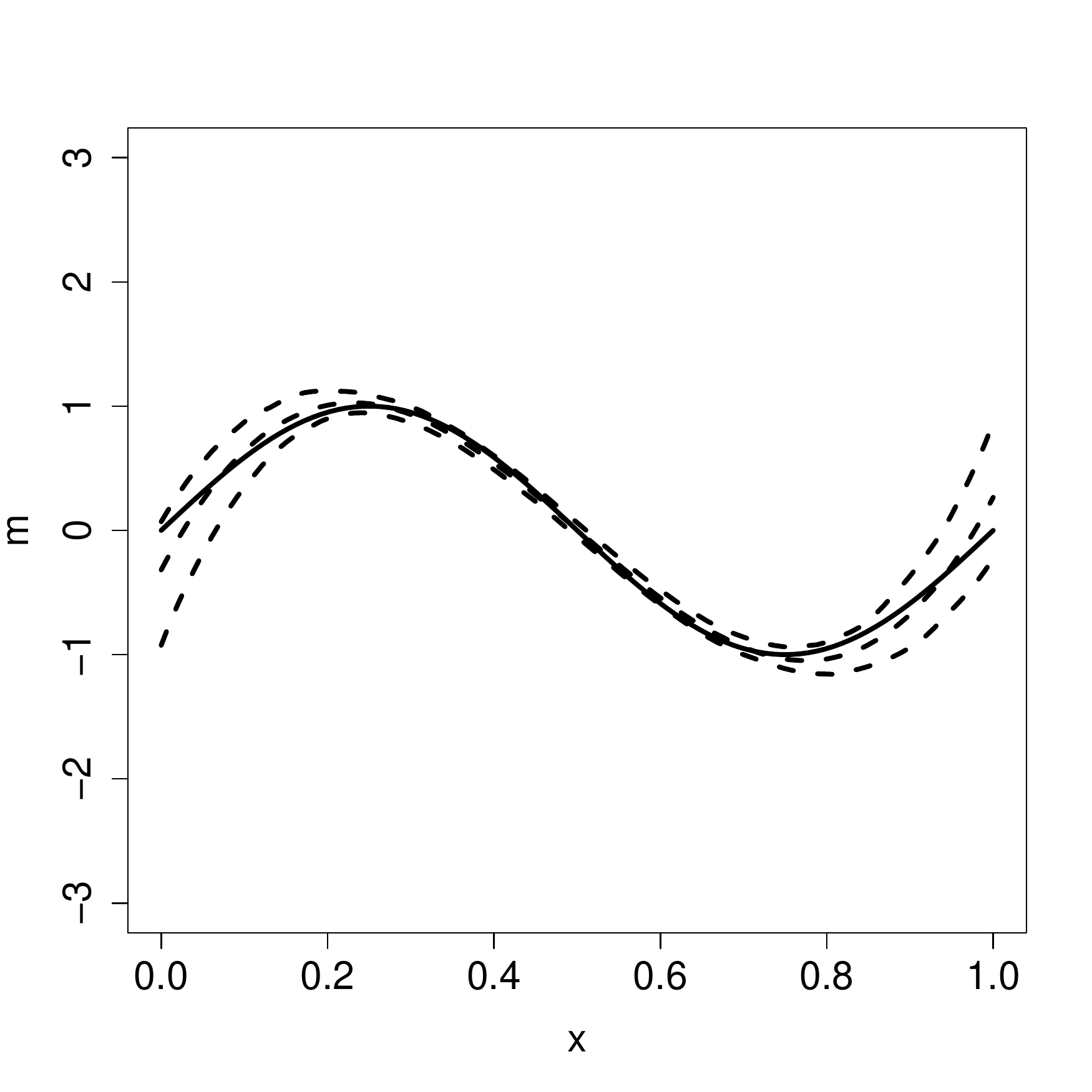}
\includegraphics[scale = 0.35]{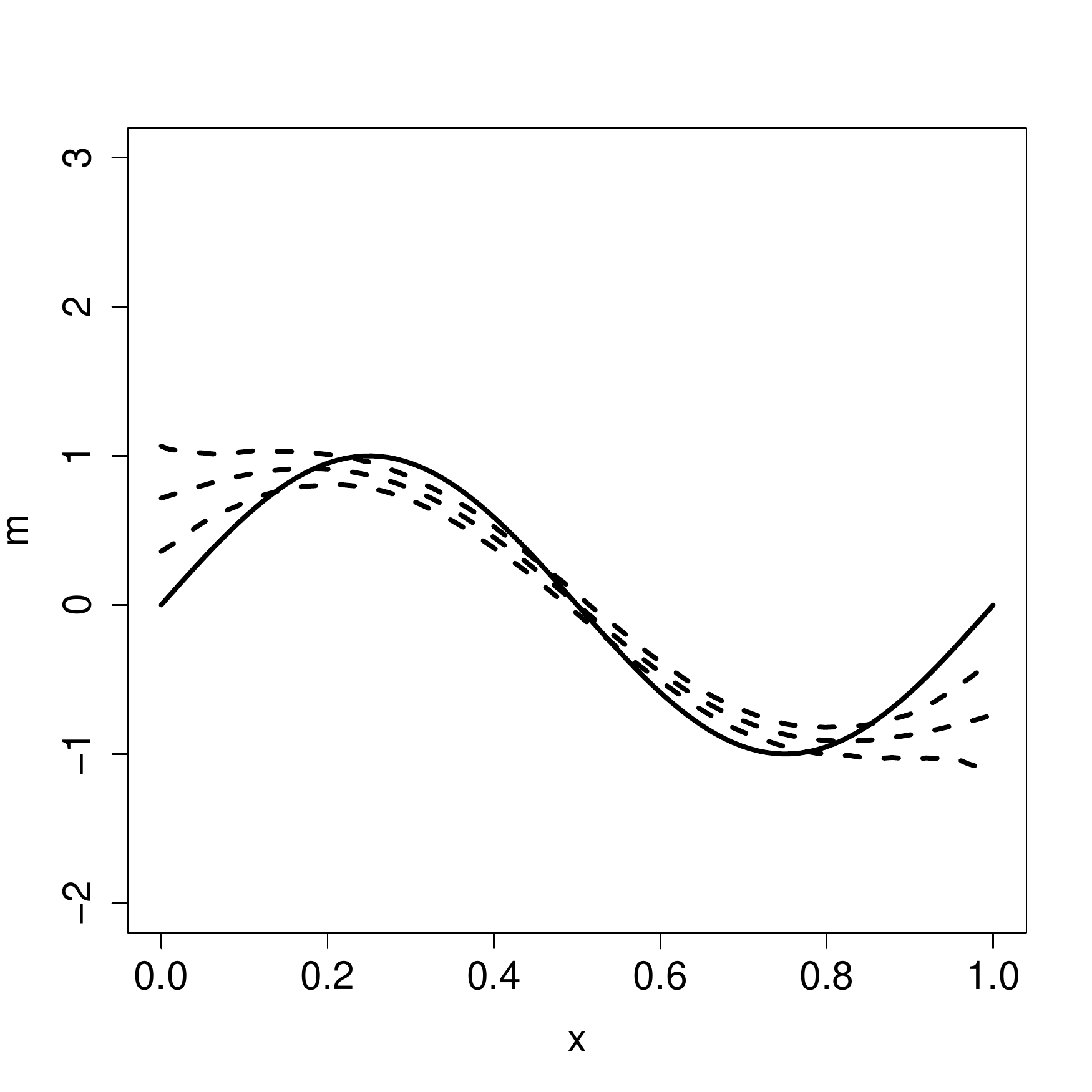}\\
\includegraphics[scale = 0.35]{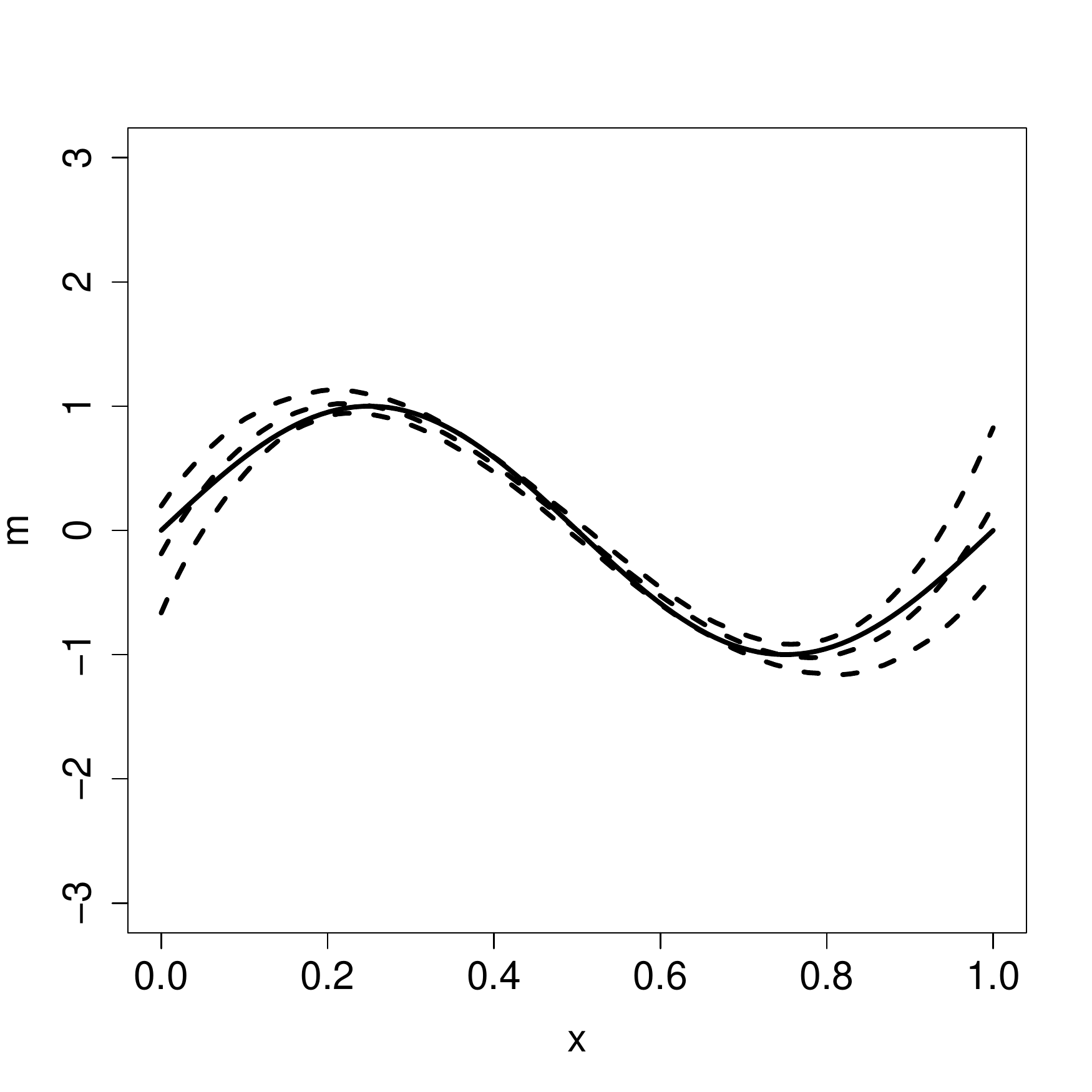}
\includegraphics[scale = 0.35]{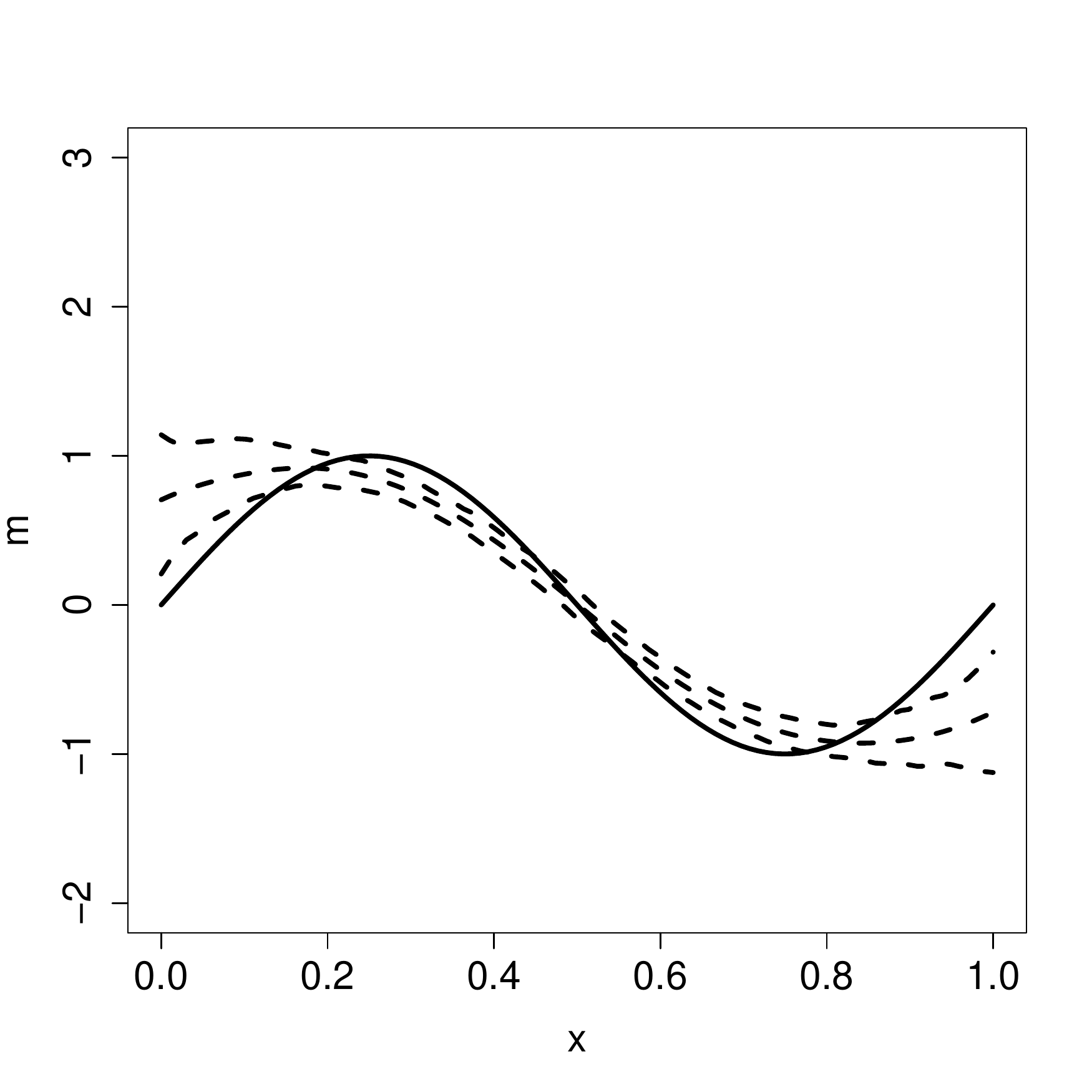}
\end{center}
\end{figure}

Similar to the pdf investigation,
we reduced the measurement
error variability and generated the errors from model II
(a)--(c) to investigate the mean function curve fitting.
We plotted the mean function estimates and the 90\%
confidence bands for the Bspline semiparametric and deconvolution
methods in Figure \ref{fig:mean1} for sample size 500, and also
provided the same results in
Figures \ref{fig:mean2}, \ref{fig:mean3} and \ref{fig:mean0} in the
Supplement.
 for sample sizes 1000, 2000, as well as 200.
These results
 indicate that the Bspline
semiparametric estimator indeed outperforms the deconvolution
method, which does not make a compact support assumption on $X$.
Their performance
difference in terms of averaged MAE between the estimated and true mean
curves are further provided in the lower half of Table
\ref{tab:compare}. The deconvolution density estimation and
nonparametric regression are implemented using the code from the website
{\em https://researchers.ms.unimelb.edu.au/$\sim$aurored/links.html.}

\begin{comment}
Similar as the pdf investigation,
we reduce the measurement
error variability and generated the errors from model II
(a)--(c) to investigate the mean function curve fitting.
We plotted the mean function estimates and the 90\%
confidence bands for the Bspline semiparametric and deconvolution
methods in Figure
\ref{fig:estmean} and Figure \ref{fig:cpreg}, respectively.
The comparison of
Figure \ref{fig:estmean} and Figure \ref{fig:cpreg} indicates that the Bspline
semiparametric estimator indeed outperforms the deconvolution method. Their performance
difference in terms of averaged MAE between the estimated and true mean
curves are further provided in the lower half of Table \ref{tab:compare}.
\end{comment}

\section{Data example}\label{sec:example}

Heavy fine particulate matter (PM2.5) air pollution has become a serious
problem in China  and its possible  effect on
respiratory diseases has been a concern in public health.
Starting from 2012, the Beijing Environmental Protection
Bureau (BEPB) has been recording the daily  PM2.5  levels in  Beijing.
Based on these data,
\cite{xu2016}
studied the effect of PM2.5 on asthma in 2013. Specifically, they explored the PM2.5 effect on the number
of daily asthma
emergency room visits (ERV) in ten hospitals in Beijing,   but found no
significant effect. In fact, the mean number
of daily asthma ERVs even shows a decreasing trend
along the increase of measured PM2.5. This
contradicts the general conclusion
that PM2.5 has short term adverse effect on asthma  \citep{fan2016}.

\begin{figure}
\centering
\caption{The estimated pdf of PM2.5 without considering measurement
error, based on data from the Beijing
Environmental Protection Bureau (solid line) and the ``Mission China'' website (dashed-line).
  }\label{fig:comden}
\includegraphics[scale = 0.3]{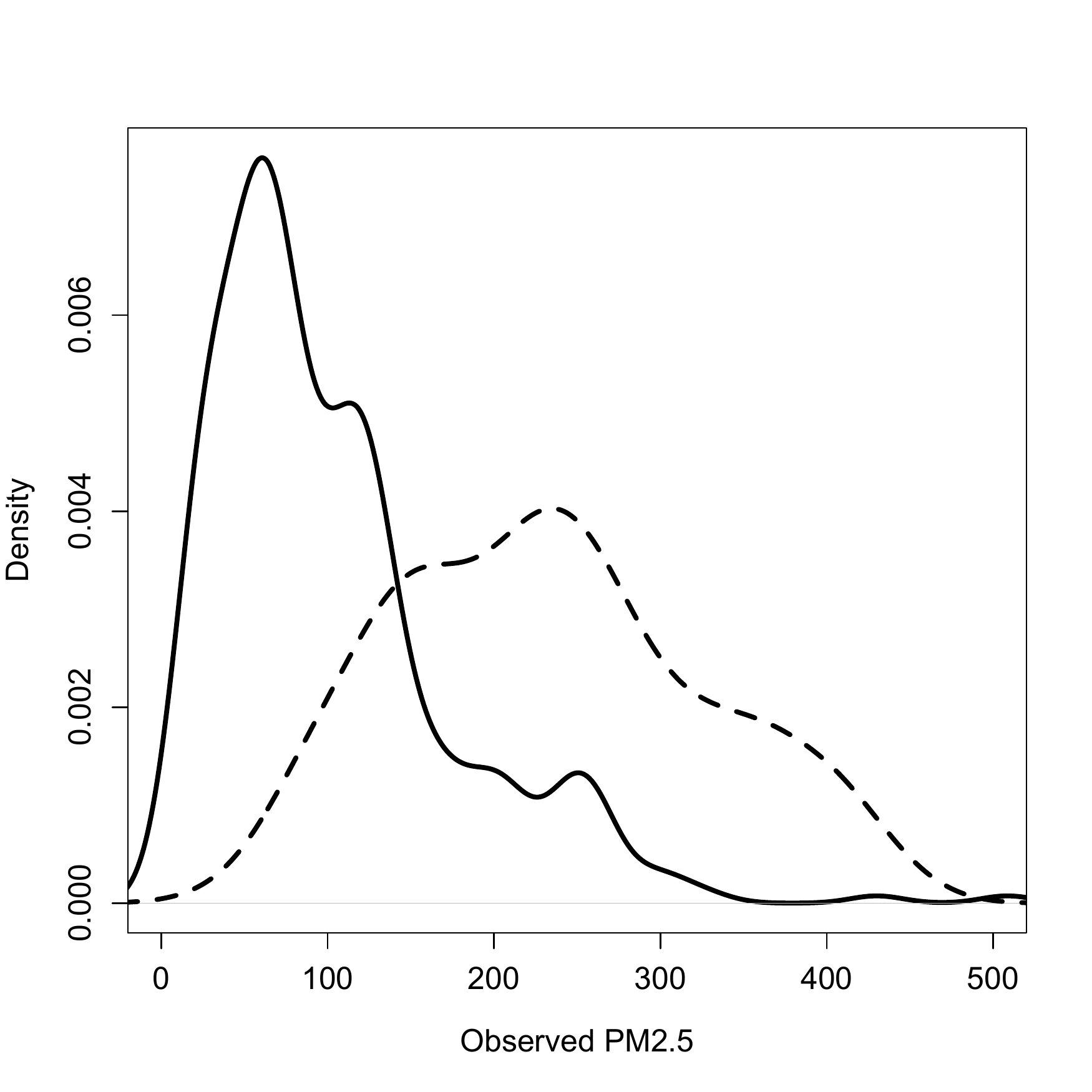}
\end{figure}

A potential reason of this inconsistency is the errors in the PM2.5
measurements which were not taken into account in the above analysis.
In fact, there are many debates on the
accuracy of the PM2.5 reports, especially in the early years
such as 2013. For example, we
compared the daily average PM2.5 reports in 2013 from  17 ambient
air quality monitoring stations and those reports from the ``Mission China
Beijing'' website \citep{WinNT} maintained by the U.S. Department of
State, and show the two estimated pdfs  of PM2.5
in Figure \ref{fig:comden}. It is clear that the estimated pdfs of
PM2.5 from the two sources are very different, where the PM2.5
concentrations obtained from the
BEPB yields  a pdf estimate with the
mode to the left of that obtained from the ``Mission China Beijing'', indicating
a generally less severe air pollution problem. This
motivates us to consider the measurement error issue in studying the effect
of PM2.5 on the daily asthma ERV.

We restrict our analysis of PM2.5
to the range from 0 to 400, since the smallest and largest recorded
PM2.5 values are
6.65 and 328.41, respectively. We first rescaled the
  observed PM2.5
by dividing $328.41-6.65$, the range of the observed PM2.5 data,  for each observation to get
scaled PM2.5. After the rescaling, the mean and variance of $\W_i$ are
0.290 and 0.045, respectively.
The data set we analyzed contains 337  observations. In the data
available to us, the $i$th observation
contains the number of
daily asthma ERVs, which we treat as response $Y_i$, the average
scaled PM2.5 level over 17 
ambient air quality monitoring stations from BEPB, which we denote as
$W_i$, and the PM2.5 level from ``Mission China Beijing'', which we
write as $W_{0i}$. To carry out the analysis, we let
 the true  scaled PM2.5 value be $X_i$, and let $W_{ki}, U_{ki}$
be the observed PM2.5 value and its corresponding measurement error at the
  $k$th monitoring station, $k = 1, \ldots, 17$. We assume $U_{ki}$'s are
  independent of each other and of $X_i$.
Then $W_i = \sum_{k =
    1}^{17}  W_{ki}/17$ and use the average $W_i$ as our surrogate
  measurement of $X_i$. To obtain the measurement error variance
  associated with $W_i$, we write $\Ubar_i =  \sum_{k =
    1}^{17}  U_{ki}/17$,  and get
$W_i= X_i + \Ubar_i$. Because our preliminary analysis result in
Figure \ref{fig:comden} suggests a possible discrepancy between the
measurements in $W_i$ and in $W_{0i}$, we allow a potential bias term
$b$ and model $W_{0i} =b +  X_i + U_{0i}$,
where $U_{0i}$ is the measurement error of $W_{0i}$. 
We assume all the
$U_{ki}$, $k=0, \dots, 17$ to have the same distribution with mean
zero, and to be independent of each other,
and we estimate $b$ by $\wh{ b} = n^{-1}(\sum_{i=1}^n W_{0i} -
\sum_{i=1}^n W_i )$, which yields the value $\wh b=0.41$.
We further estimated the variance of $\Ubar_i$ based on
$\var(\Ubar)=\{\var(W_{0i} - W_i)\}/18$. This yields
$\wh\var(\Ubar)=0.008$. Further, because $\Ubar_i$ is the average of 17
 $U_{ki}$'s, it is sensible to assume that $\Ubar_i$ has a normal
distribution.  Note that we do not assume normality on $U_{0i}$.

Based on the preliminary analysis above, we proceed to estimate
the pdf of PM2.5, i.e. $X_i$, and the mean regression function of asthma
ERVs conditional on PM2.5, i.e. $E(Y_i\mid X_i)$,
using the Bspline-assisted MLE/semiparametric methods in Sections \ref{sec:density} and
\ref{sec:regression}.
We further implemented 100 bootstraps to estimate the asymptotic variances
of the resulting estimators.
We also compared the Bspline MLE
and the Bspline semiparametric regression estimator
with the deconvolution density and regression
estimators, which does not assume the true PM2.5 range has a finite
range.  In implementing the Bspline approximation, we used two and
three equally spaced knots respectively, and in implementing the
deconvolution methods, we used bandwidth 0.05. The
number of knots is chosen based on the simulation studies in
Section \ref{sec:simulation}.

As an aside, the bandwidth 0.05 is the least we need in order to
achieve stable result for the deconvolution estimator for this data set. In fact, in selecting the
bandwidth, we implemented both the crossvalidation
\citep{stefanskicarroll1990} and the plug-in \citep{delaigle2002}
methods. Both procedures lead to very small bandwidths that induce
large numerical errors. For this reason, we increased the bandwidth to
0.05, the smallest
bandwidth that produced stable results in our
analysis.

\begin{figure}
\centering
\caption{The Bspline MLE (upper left) and deconvolution (upper right) pdf
  estimators and their 95\% confidence bands
  of PM2.5. The Bspline semiparametric (bottom left) and deconvolution (bottom right) regression
  estimators and their 95\% confidence bands.   The
    (black) solid line  are the estimators, the dashed lines are the
    95\% confidence intervals, and the (red) dot-dashed lines str the
    naive   estimator ignoring measurement errors.}\label{fig:denreg} 
\includegraphics[scale = 0.35]{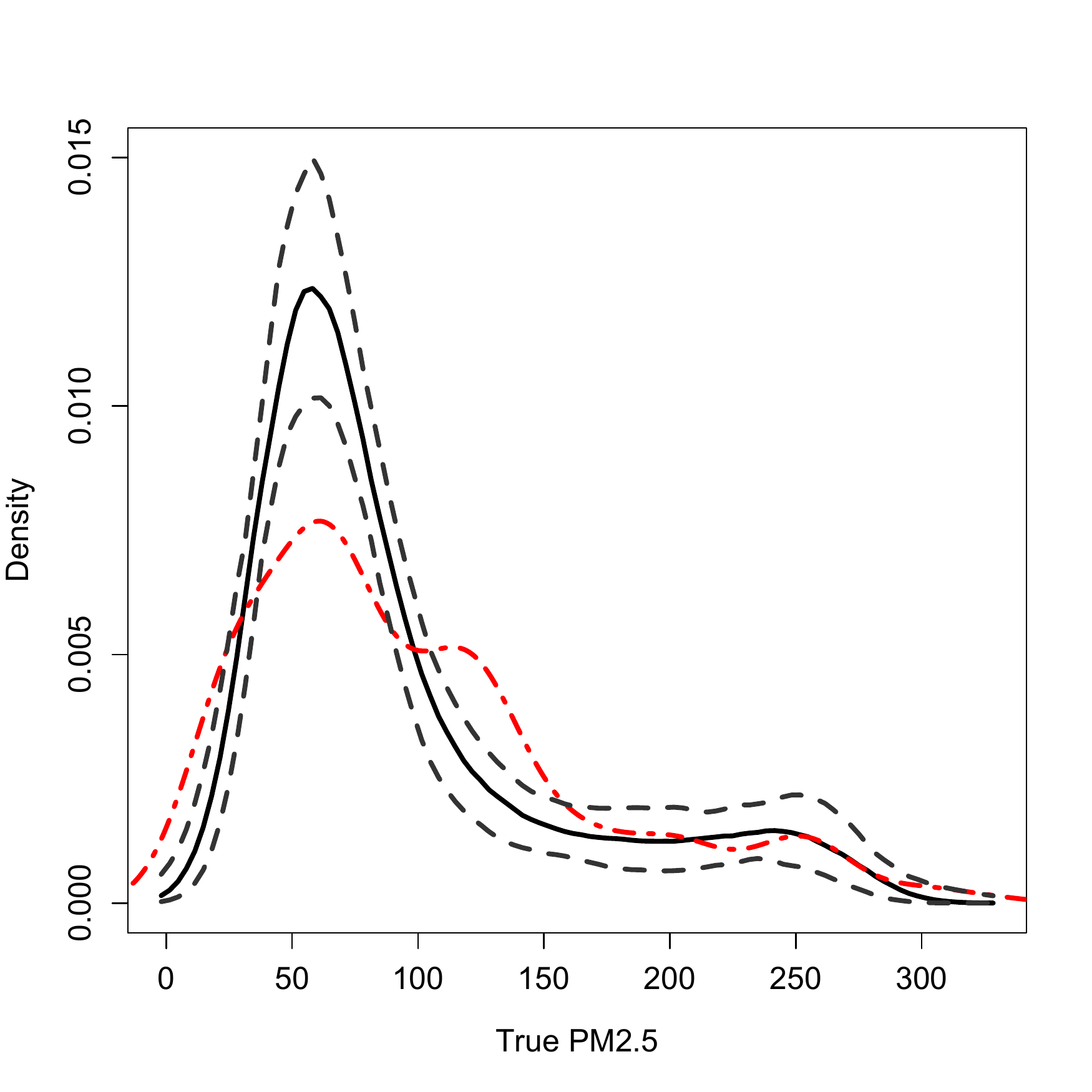}
\includegraphics[scale = 0.35]{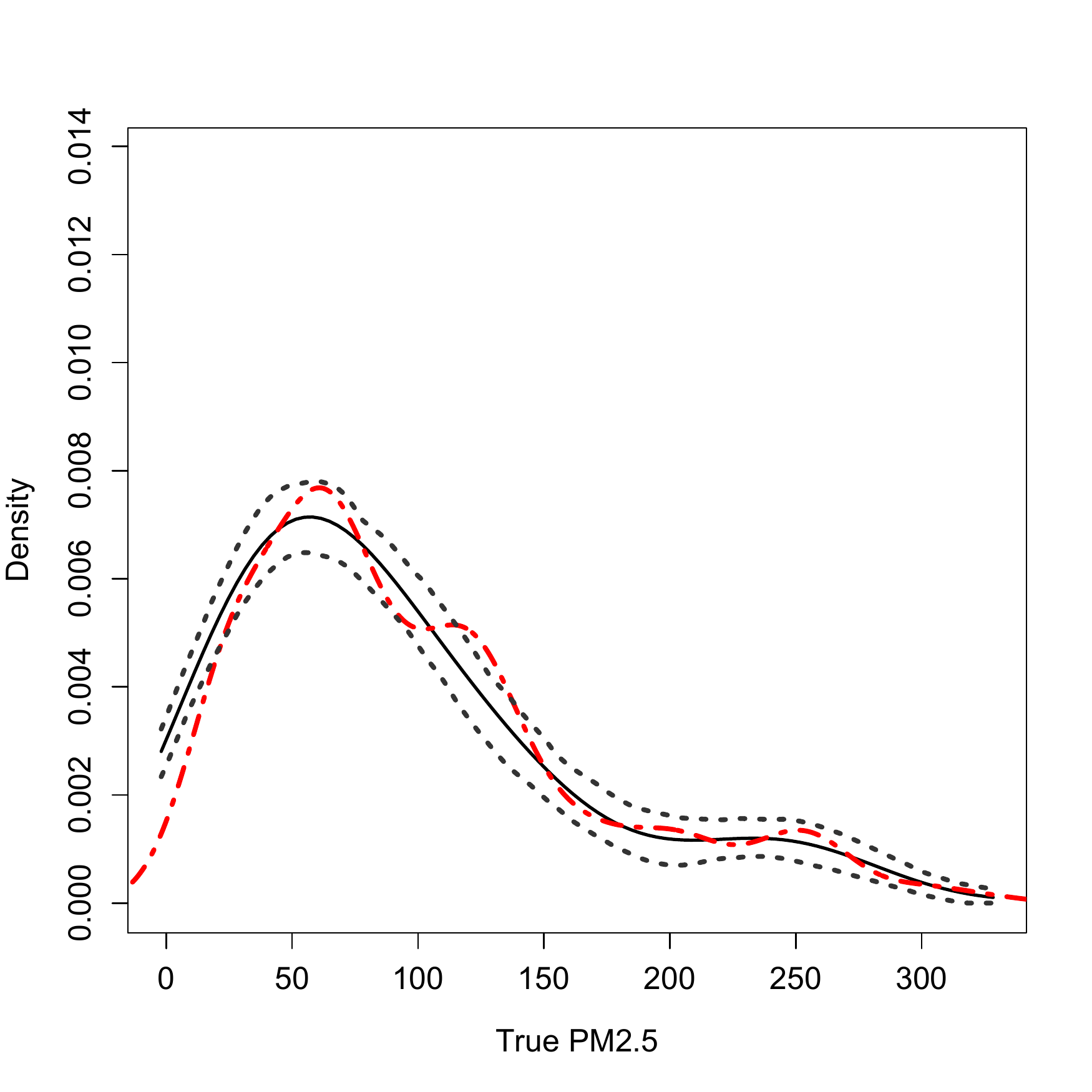}\\
\includegraphics[scale = 0.35]{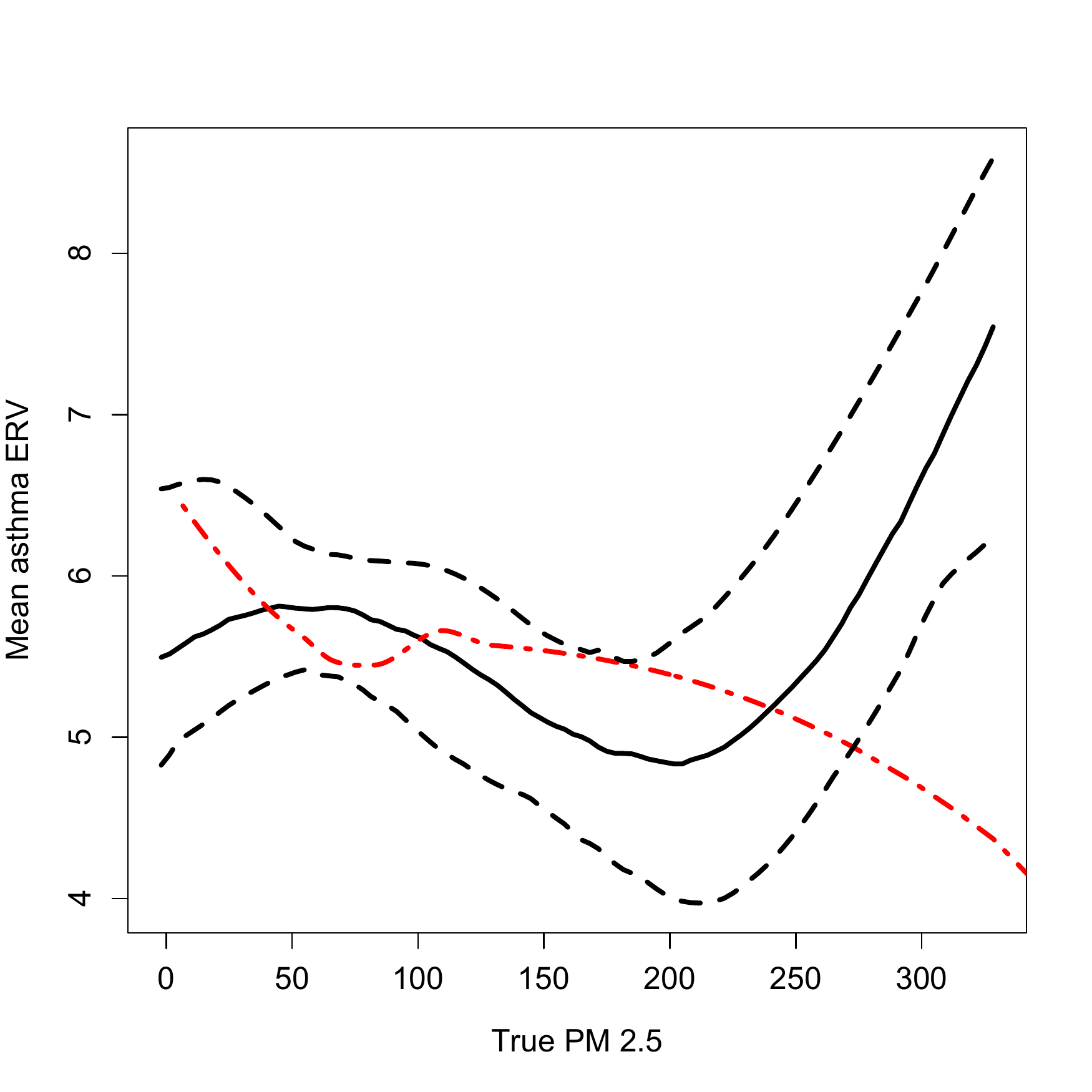}
\includegraphics[scale = 0.35]{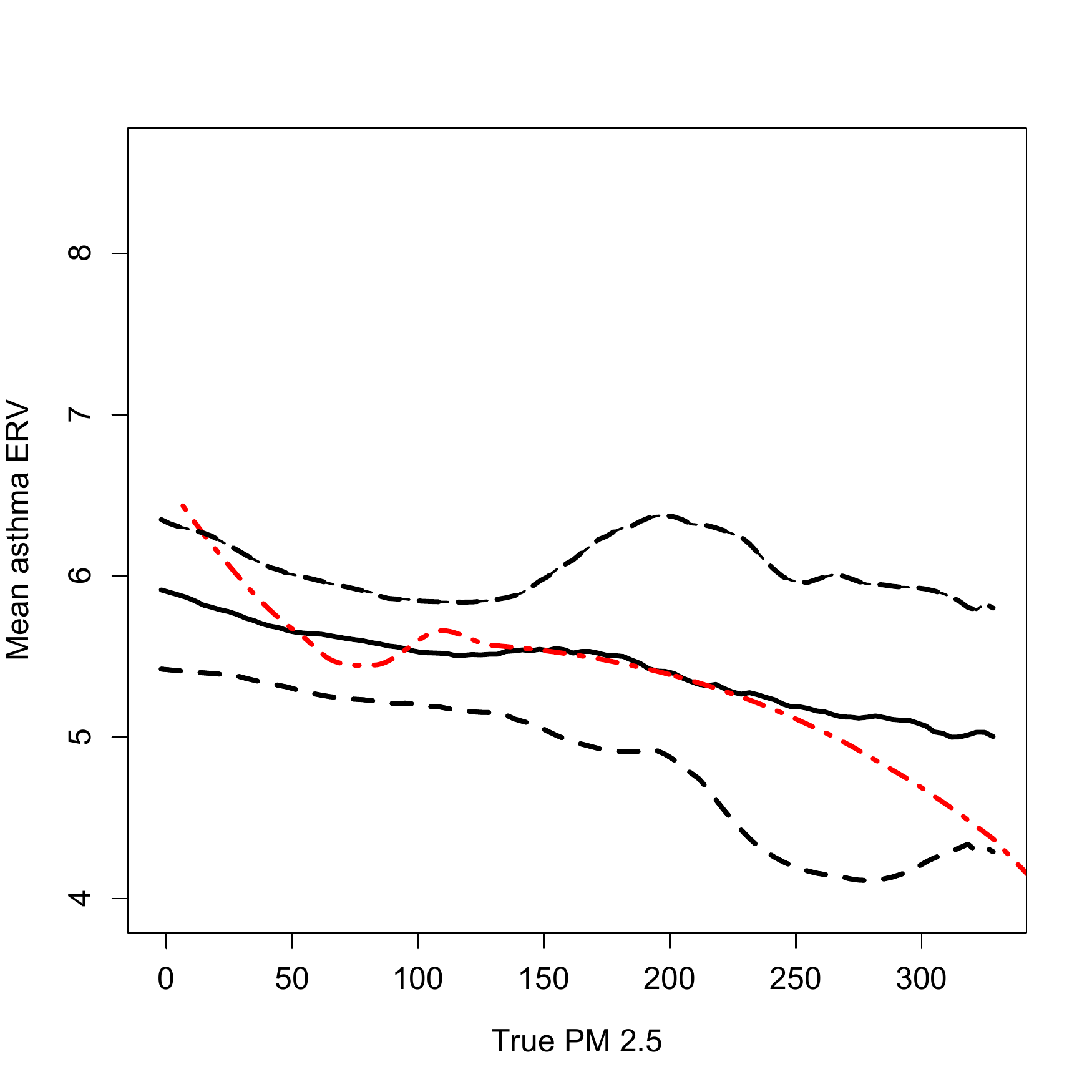}
\end{figure}

The upper panel in Figure \ref{fig:denreg} shows the estimated pdfs based on
the Bspline MLE and the deconvolution method. Compared with the kernel
estimator in the same plots which ignores the measurement errors, the
Bspline MLE shows more
difference than the deconvolution estimator. In fact, 
the noise-to-signal ratio is more than $\var(\bar
U_i)/\var(X_i)=0.25$,
hence measurement error issue is likely not ignorable.

The lower panel in
Figure \ref{fig:denreg} provides the estimated mean of
$Y_i$ as a function  of $X_i$. The Bspline semiparametric estimator
shows a fluctuating pattern in the range from 0 to 200, although
the pattern is not significant. In the range of PM2.5 concentration
larger than 200 (about 11.3\% of the observations),
it shows clearly an increasing trend,
which agrees with the conclusion in \cite{fan2016} that the
exposure to high PM2.5 has an adverse effect on  asthma
onset. In contrast, the relation from the deconvolution estimator is
similar to that of the
local linear regression estimator ignoring the measurement errors,
and it is unable to detect the increasing trend of the asthma ERVs
as the PM2.5 level increases.

\section{Discussion}\label{sec:discussion}

We have developed a Bspline-assisted MLE for nonparametric pdf estimation and a
Bspline-assisted semiparametric estimator for nonparametric
regression mean function estimation,
in the situation that the covariates are compactly supported and
measured with error. The
performance of both procedures are superior to the widely used
deconvolution methods that do not assume compact support on the
covariates, in terms of both their theoretical convergence
rate and their numerical performance.

A key difference between our procedures and the deconvolution
procedures is that we restrict our interest in estimating the functions on a
compact set while the deconvolution method does not impose such an
assumption. Practically, the
possible range of a random variable is indeed finite,  and the relevant
information is needed only for functions in a finite range, we consider
the compact support assumption to be mild.
We are thus very curious if deconvolution
methods can achieve the same convergence property in such case,
with possibly some modifications on the existing procedures. To this
end, \cite{hall2005} provides some relevant results based on a discrete
Fourier transform and its inverse, while
we leave this general question as an open problem for researchers who specialize
in deconvolution methods.

Our method is the first non-deconvolution procedure that
achieves the standard nonparametric convergence rate, the same rate as that obtained
without measurement error. It would be interesting to investigate if other methods,
such as kernel or Fourier series methods, can also attain such rate.

The density and regression function estimation problems
  studied in this work are the most basic ones. In practice, various
  complications may occur. For example, the error distribution
  $f_U(u)$ may not be known and need to be estimated parametrically or
  nonparametrically based on repeated measurements or other
  instruments. If one estimates $f_U(\cdot)$ parametrically first then insert it into
  our procedure, it will not have first order effect. However, if
  $f_U(\cdot)$ is estimated nonparametrically, it will in general
  affect both the bias and  variance of the resulting estimation of
  both the density and regression functions. See
\cite{delaigle2016} for
deconvolution based estimation incorporating an unknown error distribution.

%\pagebreak
%\newpage
%\setcounter{page}{1}
\section*{Appendix}
\renewcommand{\theequation}{A.\arabic{equation}}
\renewcommand{\thesubsection}{A.\arabic{subsection}}
\setcounter{equation}{0}
\setcounter{subsection}{0}

\subsection{Proof of Theorem \ref{th:density}}\label{sec:appdensity}
Recall that, $\a_L$  is the
subvector of $\a$ without the first element.

It follows that
\bse
&&\sup_{x \in [0, 1]} |\log\{\wh f_X(x)\}
- \log \{f_{X0}(x)\}|\\
&=&\sup_{x \in [0, 1]} |\B\trans_r(x)\wh{\btheta}
-\log\int_0^1\exp\{\B\trans_r(x) \wh{\btheta}\}dx
- \log \{f_{X0}(x)\}|\\
&=& \sup_{x \in [0, 1]} |\B\trans_r(x)\wh{\btheta} - \B\trans_r(x){\btheta}_0
+ \B\trans_r(x){\btheta}_0
- \log \{f_{X0}(x)\}\\
&&
-\log\int_0^1\exp\{\B\trans_r(x) \wh{\btheta}\}dx
+\log\int_0^1\exp\{\B\trans_r(x) \btheta_0\}dx  \\
&&-\log\int_0^1\exp\{\B\trans_r(x) \btheta_0\}dx|
\ese
\bse
&\leq & \sup_{x \in [0, 1]} \{ \|\B_r(x)\|_2  \|\wh{\btheta} - {\btheta}_0\|_2\}
 +\sup_{x \in [0, 1]} |\B\trans_r(x){\btheta}_0 - \log\int_0^1\exp\{\B\trans_r(x) {\btheta}_0\}dx-
 \log\{f_{X0}(x)\}|\\
&&
 + |\log\int_0^1\exp\{\B\trans_r(x) \wh{\btheta}\}dx -
 \log\int_0^1\exp\{\B\trans_r(x) {\btheta}_0\}dx |\\
&=& \sup_{x \in [0, 1]} \{ \|\B_r(x)\|_2  \|\wh{\btheta} - {\btheta}_0\|_2\}
 +\sup_{x \in [0, 1]} |\B\trans_r(x){\btheta}_0
 -\log\int_0^1\exp\{\B\trans_r(x) {\btheta}_0\}dx-
 \log\{f_{X0}(x)\}|\\
&&
 + |\frac{\partial}{\partial \btheta\trans} \log\int_0^1\exp\{\B\trans_r(x)
 {\btheta}^\star\}dx (\wh{\btheta} - \btheta_0) | \\
&\leq& \sup_{x \in [0, 1]} \{ \|\B_r(x)\|_2  \|\wh{\btheta} - {\btheta}_0\|_2\}
 +\sup_{x \in [0, 1]} |\B\trans_r(x){\btheta}_0
 -\log\int_0^1\exp\{\B\trans_r(x) {\btheta}_0\}dx-
 \log\{f_{X0}(x)\}|\\
&&
 + \|\frac{\partial}{\partial \btheta\trans} \log\int_0^1\exp\{\B\trans_r(x)
 {\btheta}^\star\}dx \|_2\|(\wh{\btheta} - \btheta_0) \|_2 \\
&=& O_p\{(nh_b)^{-1/2} + h_b^q\},
\ese
where $\btheta^\star$ is the point on the line connecting
$\wh{\btheta}$ and $\btheta_0$.
The last equality holds by Proposition \ref{pro:theta} and Condition \ref{ass:dis} that
\bse
&&|\B\trans_r(x){\btheta}_0
 -\log\int_0^1\exp\{\B\trans_r(x) {\btheta}_0\}dx-
 \log\{f_{X0}(x)\}| = O_p(h_b^q).
\ese
This also implies $\sup_{x \in [0, 1]} | \wh f_X(x)
- f_{X0}(x)| = O_p\{(nh_b)^{-1/2} + h_b^q\}$.

In addition
 \bse
&&\wh{f}_X(x)
-f_{X0}(x) \\
&= &\left[\frac{\exp\{\B\trans_r(x)\wh{\btheta}\}}
{\int_0^1\exp\{\B\trans_r(x) \wh{\btheta}\}dx}
-f_{X0}(x)\right]\\
&=&\left[\frac{\exp\{\B\trans_r(x)\wh{\btheta}\}}
{\int_0^1\exp\{\B\trans_r(x) \wh{\btheta}\}dx} - \frac{\exp\{\B_r(x)
  \trans{\btheta}_0\}}
{\int_0^1\exp\{\B\trans_r(x) {\btheta}_0\}dx}\right.\\
&&\left.+  \frac{\exp\{\B\trans_r(x){\btheta}_0\}}
{\int_0^1\exp\{\B\trans_r(x) {\btheta}_0\}dx}
-f_{X0}(x)\right]\\
&=&\left[\frac{\exp\{\B\trans_r(x)\wh{\btheta}\}}
{\int_0^1\exp\{\B\trans_r(x) \wh{\btheta}\}dx} - \frac{\exp\{\B_r(x)
  \trans{\btheta}_0\}}
{\int_0^1\exp\{\B\trans_r(x) {\btheta}_0\}dx}\right]+ O_p(h_b^q)\\
&=&\frac{\partial}{\partial\btheta_L\trans} \left[\frac{\exp\{\B\trans_r(x)\btheta_0\}}
{\int_0^1\exp\{\B\trans_r(x) \btheta_0 \}dx}\right]
\left(E \left[\frac{\int_0^1 f_{X0}(x)f_U(W_i-x) \B_{rL}(x)\B\trans_{rL}(x)
dx }{\int_0^1 f_{X0}(x)f_U(W_i-x) dx}\right.\right.\\
&&- \frac{\left\{\int_0^1f_{X0}(x)f_U(W_i-x) \B_{rL}(x)
dx \right\}^{\otimes2}}{\left\{\int_0^1 f_{X0}(x)f_U(W_i-x) dx \right\}^2}- \int_0^1 f_{X0}(x)\B_{rL}(x)
\B\trans_{rL}(x) dx\\
&&\left.\left.+ \left\{\int_0^1 f_{X0}(x)\B_{rL}(x)
  dx\right\}^{\otimes2}\right]\right)^{-1}  n^{-1}\sumi
\int_0^1\left\{\frac{f_{X0}(x)f_U(W_i-x)
  }{\int_0^1f_{X0}(x)f_U(W_i-x) dx}\right.\\
&&\left.- f_{X0}(x)\right\}\B_{rL}(x)dx\{1 + o_p(1)\}+ O_p(h_b^q).
\ese
The last equality holds by Proposition \ref{pro:theta}. 
Taking expectations on both sides, because the above leading term has
mean 0 and is of order $O_p\{(nh_b)^{-1/2}\}$, we have
\bse
\bias\{\wh{f}_X(x)\} = O_p(h_b^q) + o_p\{(nh_b)^{-1/2}\}.
\ese
Therefore we have
\bse
\sqrt{nh_b} [\wh{f}_X(x)
-f_{X0}(x) -\bias\{\wh{f}_X(x)\} ]= n^{-1/2}\sumi
Q(W_i,x)+o_p(1).
\ese
and
when $Nn^{-1/(2q + 1)} \to \infty$ we further have
 \bse
\sqrt{nh_b} \{\wh{f}_X(x)
-f_{X0}(x) \} = n^{-1/2}\sumi
Q(W_i,x)+o_p(1).
\ese
This proves Theorem \ref{th:density}. \qed

\subsection{Proof of Theorem \ref{th:mean}}\label{sec:appmean}
\bse
&&\sup_{x \in [0, 1]} |\wh{m}(x) - m(x)| \\
&=& \sup_{x \in [0, 1]} |\B\trans_r(x)\wh{\bb} - \B\trans_r(x) {\bb_0}
+\B\trans_r(x) {\bb_0} -   m(x)|\\
&\leq& \sup_{x \in [0, 1]} |\B\trans_r(x)\wh{\bb} - \B\trans_r(x) {\bb_0}|
+\sup_{x \in [0, 1]} |\B\trans_r(x) {\bb_0} -   m(x)|\\
&=& O_p\{(nh_b)^{-1/2} + h_b^q\},
\ese
by Proposition \ref{pro:theta} and Condition \ref{ass:mdis}.

Further, we have
\bse
&&\wh{m}(x) - m(x) \\
&=& \B\trans_r(x)\\
&&\times\left[ - \left\{E\left(\frac{\partial[\S^*_{\bb}(W_i, Y_i, \bb)
    -E^*\{\a(X, \bb) \mid Y_i, W_i,
\bb\}]}{\partial \bb\trans}\bigg |_{\B_r(\cdot)\trans\bb =
m(\cdot)}\right)\right\}^{-1}\right.\\
&&\left.\times  \frac{1}{n}\sumi \S^*_{\bb}(W_i, Y_i, m) -E^*\{\a(X, m) \mid Y_i, W_i,
m \} \right] \\
&= & \B\trans_r(x)\\
&&\times\left[ - \left\{E\left(\frac{\partial[\S^*_{\bb}(W_i, Y_i, \bb)
    -E^*\{\a(X, \bb) \mid Y_i, W_i,
\bb\}]}{\partial \bb\trans}\bigg |_{\B_r(\cdot)\trans\bb =
m(\cdot)}\right)\right\}^{-1}\right.\\
&&\left.\times  n^{-1}\sumi \S^*_{\bb}(W_i, Y_i, m) -E^*\{\a(X, m) \mid Y_i, W_i,
m \} \right] \{1 + o_p(1)\}+ O_p(h_b^q). \\
\ese
The last equality holds by Proposition \ref{pro:beta}.
Taking expectations on both sides, because the above leading term has
mean 0 and is of order $O_p\{(nh_b)^{-1/2}\}$,  we have
\bse
\bias\{\wh{m}(x) \} = O_p(h_b^q) + o_p\{(nh_b)^{-1/2}\}.
\ese
Therefore,
\bse
\sqrt{nh_b}[\wh{m}(x) - m(x) - \bias\{\wh{m}(x)\} ]= n^{-1/2}\sumi
Q(W_i,Y_i,x)+o_p(1).
\ese
When $N n^{1/(2q + 1)} \to
\infty$, we further have
\bse
\sqrt{nh_b}\{\wh{m}(x) - m(x) \} = n^{-1/2}\sumi
Q(W_i,Y_i,x)+o_p(1).
\ese
This proves the result.
\qed

\section*{Acknowledgments}

We thank Yang Li for providing the daily asthma emergency room
visits data. Ma's research was partially supported by grants from NSF (1608540)
and NIH (HL138306, NS073671).
  Carroll's research was supported by a grant from the National Cancer
  Institute (U01-CA057030).

%\newpage

\clearpage
\begin{table}[!h]
\centering
\caption{Performance comparison between the Bspline
  MLE/semiparametric estimator and deconvolution. Mean of the maximum absolute
  differences are reported.
Results based
on average over 200 simulations.
}\label{tab:compare}
\begin{tabular}{ccccccc}
\hline
\multicolumn{7}{c}{\hskip 5mm}\\
\multicolumn{7}{c}{pdf estimation: $ E \{\sup_x |\wh f_{X}(x)-f_{X0}(x)| \}$}\\
&\multicolumn{3}{c}{Bspline MLE}&\multicolumn{3}{c}{Deconvolution}\\
&$n = 500$& $n = 1000$&$n = 2000$&$n = 500$&$n = 1000$&$n = 2000$\\
Model II(a)&0.203 &0.156& 0.107&0.238& 0.212& 0.160\\
Model II(b)&0.213& 0.155& 0.104&0.230& 0.197& 0.158\\
Model II(c)&0.230 &0.181& 0.131&0.315& 0.242& 0.230\\
\hline
\multicolumn{7}{c}{\hskip 5mm}\\
\multicolumn{7}{c}{mean estimation: $ E\{\sup_x |\wh m(x)-m(x)| \}$}\\
&\multicolumn{3}{c}{Bspline semiparametric}&\multicolumn{3}{c}{Deconvolution}\\
&$n = 500$& $n = 1000$&$n = 2000$&$n = 500$&$n = 1000$&$n = 2000$\\
Model II(a)& 0.370 &0.263 &0.175&0.908& 0.796 &0.762\\
Model II(b)&0.425 &0.264 &0.163&0.857 &0.828& 0.779\\
Model II(c)&0.414 &0.291 &0.219&0.880& 0.832& 0.801\\
\hline
\end{tabular}
\end{table}

%\clearpage

\newpage

\bibliographystyle{agsm}
\bibliography{errorrate}

\end{document}